\documentclass[12pt, leqno]{article}

\usepackage{amsmath, amssymb, amsthm}
\usepackage{fullpage}

\newtheorem{theorem}{Theorem}
\newtheorem{corollary}{Corollary}[theorem]
\newtheorem{lemma}{Lemma}[section]

\newtheorem*{Goldbach}{Goldbach's Conjecture}
\newtheorem*{GRH}{Generalized Riemann Hypothesis}
\newtheorem*{RH}{Riemann Hypothesis}

\theoremstyle{remark}
\newtheorem*{remark}{Remark}
\newtheorem*{conclusion}{Conclusion}
\newtheorem*{acknowledgement}{Acknowledgements}

\numberwithin{equation}{section}

\newcommand{\bysame}{\rule{.3in}{.5pt}}
\newcommand{\eps}{\varepsilon}
\newcommand{\modulo}[1]{\mathrm{mod}\;#1}
\newcommand{\pmodulo}[1]{\;(\mathrm{mod}\;#1)}

\DeclareMathOperator{\re}{Re}
\DeclareMathOperator{\li}{li}

\title{An invitation to additive prime number theory}
\author{A. V. Kumchev and D. I. Tolev}
\date{\today}

\begin{document}

\maketitle 

\begin{abstract}
  The main purpose of this survey is to introduce the inexperienced reader to
additive prime number theory and some related branches of analytic number theory.
We state the main problems in the field, sketch their history and the basic
machinery used to study them, and try to give a representative sample of the
directions of current research.

\end{abstract}

\section{Introduction}
\label{s1}

Additive number theory is the branch of number theory that studies
the representations of natural numbers as sums of integers subject
to various arithmetic restrictions. For example, given a sequence
of integers 
\[
  \mathcal A = \{ a_1 < a_2 < a_3 < \cdots \}
\]
one often asks
what natural numbers can be represented as sums of a fixed number
of elements of $\mathcal A$; that is, for any fixed $s \in \mathbb
N$, one wants to find the natural numbers $n$ such that the
diophantine equation
\begin{equation} \label{1.1}
  x_1 + \dots + x_s = n
\end{equation}
has a solution in $x_1, \dots, x_s \in \mathcal A$. The sequence
$\mathcal A$ may be described in some generality (say, one may
assume that $\mathcal A$ contains ``many" integers), or it may
be a particular sequence of some arithmetic interest (say, $\mathcal
A$ may be the sequence of $k$th powers, the sequence of prime
numbers, the values taken by a polynomial $F(X) \in \mathbb Z[X]$
at the positive integers or at the primes, etc.). In this survey,
we discuss almost exclusively problems of the latter kind.
The main focus will be on two questions, known as Goldbach's
problem and the Waring--Goldbach problem, which are concerned
with representations as sums of primes and powers of primes,
respectively.

\subsection{Goldbach's problem}
\label{s1.1}

Goldbach's problem appeared for the first time in 1742 in the
correspondence between Goldbach and Euler. In modern language,
it can be stated as follows.

\begin{Goldbach}
  Every even integer $n \ge 4$ is the sum of two primes, and 
  every odd integer $n \ge 7$ is the sum of three primes.
\end{Goldbach} 

The two parts of this conjecture are known as the binary Goldbach
problem and the ternary Goldbach problem, respectively. Clearly,
the binary conjecture is the stronger one. It is also much more
difficult.
 
The first theoretical evidence in support of Goldbach's conjecture
was obtained by Brun~\cite{Brun19}, who showed that every large even 
integer is the sum of two integers having at most nine prime factors. 
Brun also obtained an upper bound of the correct order for the number 
of representations of a large even integer as the sum of two primes.

During the early 1920s Hardy and Littlewood~\cite{HaLi19}--%
\cite{HaLi25} developed the ideas in an earlier paper by Hardy and 
Ramanujan~\cite{HaRa18} into a new analytic method in additive 
number theory. Their method is known as the {\em circle method}. In 
1923 Hardy and Littlewood \cite{HaLi23a, HaLi23b} applied the circle 
method to Goldbach's problem. Assuming the Generalized Riemann 
Hypothesis\footnote{An important conjecture about certain 
Dirichlet series; 
see \S \ref{ssiii2} for details.} (GRH), they proved that all but 
finitely many odd integers are sums of three primes and that all but
$O \left( x^{1/2 + \eps} \right)$ even integers $n \le x$ are sums of 
two primes. (Henceforth, $\eps$ denotes a positive number which can 
be chosen arbitrarily small if the implied constant is allowed to 
depend on $\eps$.)

During the 1930s Schnirelmann~\cite{Schn33} developed 
a probabilistic approach towards problems in additive number theory.
Using his method and Brun's results, he was able to prove 
unconditionally that there exists a positive integer $s$ such that 
every sufficiently large integer is the sum of at most $s$ primes. 
Although the value of $s$ arising from this approach is much larger 
than the conjectured $s = 3$, Schnirelmann's result represented a 
significant achievement, as it defeated the popular belief at the
time that the solution of Goldbach's problem must depend on GRH.
(Since its first appearance, Schnirelmann's method has been polished
significantly. In particular, the best result to date obtained in
this fashion by Ramare \cite{Rama95} states that one can take $s=7$.)

In 1937 I.~M.~Vinogradov~\cite{IVin37a} found an ingenious new
method for estimating sums over primes, which he applied to the
exponential sum
\begin{equation}\label{1.2}
  f(\alpha) = \sum_{p \le n} e( \alpha p ),
\end{equation}
where $\alpha$ is real, $p$ denotes a prime, and $e(\alpha) =
\exp \left( 2\pi {\mathrm i} \alpha \right)$. Using his estimate 
for $f(\alpha)$, Vinogradov was able to give a new, unconditional 
proof of the result of Hardy and Littlewood on the ternary Goldbach
problem. His result is known as Vinogradov's three prime theorem.

\begin{theorem}[Vinogradov, 1937]\label{th1}
  For a positive integer $n$, let $R(n)$ denote the number of 
  representations of $n$ as the sum of three primes. Then
  \begin{equation} \label{1.3}
    R(n) = \frac {n^2}{2(\log n)^3} \mathfrak S(n)
    + O \left( n^2 (\log n)^{-4} \right),
  \end{equation}
  where
  \begin{equation} \label{1.4}
    \mathfrak S(n) = 
    \prod_{p \mid n} \left( 1 - \frac {1}{(p-1)^2} \right)
    \prod_{p \nmid n} \left( 1 + \frac {1}{(p-1)^3} \right).
  \end{equation}
  In particular, every sufficiently large odd integer is the sum
  of three primes.
\end{theorem}

The products in \eqref{1.4} are over the primes dividing $n$ and
over those not dividing $n$, respectively. In particular, when $n$
is even, we have $\mathfrak S(n) = 0$, making \eqref{1.3} trivial.
On the other hand, when $n$ is odd, we have $\mathfrak S(n) \ge 
1$. We describe the proof of Theorem~\ref{th1} in \S \ref{ssii1}.

It should be noted that the independence of GRH in Theorem \ref{th1}
comes at the price of a mind-boggling implied constant. If one 
avoids $O$-notation and makes all the constants explicit, one finds
that the original (GRH-dependent) work of Hardy and Littlewood
establishes the ternary Goldbach conjecture for $n \ge 10^{50}$,
whereas Vinogradov's method requires $n \ge 10^{6\,800\,000}$ and even 
its most refined version available today (see Liu and Wang
\cite{LiWa02}) requires $n \ge 10^{1\,346}$. To put these numbers in 
perspective, we remark that even the bound $10^{50}$ is beyond hope 
of ``checking the remaining cases by a computer". In fact, only 
recently have Deshouillers {\em et al.}~\cite{DeEfRiZi97} proved that 
if GRH is true, the ternary Goldbach conjecture holds for all odd 
$n \ge 7$.

In 1938, using Vinogradov's method, Chudakov~\cite{Chud38}, van der 
Corput \cite{vdCo38}, and Estermann~\cite{Este38} each showed 
that {\em almost all} even integers $n \le x$ are sums of two primes. 
More precisely, they proved  that for any $A > 0$ we have
\begin{equation}\label{1.5}
  E(x) = O \left( x(\log x)^{-A} \right),
\end{equation}
where $E(x)$ denotes the number of even integers $n \le x$ that cannot 
be represented as the sum of two primes. The first improvement on 
\eqref{1.5} was obtained by Vaughan~\cite{Vaug72}.  It was followed by
a celebrated work by Montgomery and Vaughan~\cite{MoVa75} from 1975, 
in which they established the existence of an absolute constant
$\delta > 0$ such that
\begin{equation}\label{1.6}
  E(x) = O \left( x^{1 - \delta} \right).
\end{equation}
The first to compute an explicit numerical value for $\delta$
were Chen and Pan~\cite{ChPa80}. They showed that the method of 
Montgomery and Vaughan yields \eqref{1.6} with $\delta = 0.01$.
Subsequently, this result has been sharpened by several authors
and currently \eqref{1.6} is known to hold with $\delta = 0.086$
(see Li~\cite{HLi00}). In June 2004, Pintz \cite{Pint04} 
announced a further improvement on \eqref{1.6}. He has 
established the above bound with $\delta = \frac 13$ and can 
also show that for all but $O(x^{3/5 + \eps})$ even integers 
$n \le x$ either $n$ or $n - 2$ is the sum of two primes.

One may also think of the binary Goldbach conjecture as a claim
about the primes in the sequence
\begin{equation}\label{1.7}
  \mathcal A = \mathcal A(n)
  = \left\{ n - p : p \text{ prime number}, \; 2 < p < n \right\},
\end{equation}
namely, that such primes exist for all even $n \ge 6$. Denote by
$P_r$ an integer having at most $r$ prime factors, counted with
their multiplicities, and refer to such a number as an {\em almost
prime} of order $r$ (thus, Brun's result mentioned above asserts 
that every large even $n$ can be represented in the form $n = P_9 
+ P_9'$). In 1947 R\'enyi~\cite{Reny47} proved that there is a 
fixed integer $r$ such that the sequence $\mathcal A$ contains 
a $P_r$-number when $n$ is sufficiently large. Subsequent work by 
many mathematicians reduced the value of $r$ in R\'enyi's result 
almost to the possible limit and fell just short of proving the 
binary Goldbach conjecture. The best result to date was obtained 
by Chen~\cite{Chen73}.

\begin{theorem}[Chen, 1973]\label{th2}
  For an even integer $n$, let $r(n)$ denote the number of 
  representations of $n$ in the form $n = p + P_2$, where $p$ is
  a prime and $P_2$ is an almost prime of order $2$. There exists
  an absolute constant $n_0$ such that if $n \ge n_0$, then
  \[
    r(n) > 0.67 \prod_{p > 2} \left( 1 - \frac 1{(p - 1)^2} \right)
    \prod_{ \substack{ p > 2\\ p \mid n}}
    \left( \frac {p - 1}{p - 2} \right) \frac n{(\log n)^2}.
  \]
  In particular, every sufficiently large even integer $n$ can be
  represented in the form $n = p + P_2$.
\end{theorem}

\subsection{Waring's problem}
\label{s1.2}

Before proceeding with the Waring--Goldbach problem, we will
make a detour to present the most important results in Waring's
problem, as those results and the work on Goldbach's problem
have been the main motivation behind the Waring--Golbach problem.
It was probably the ancient Greeks who first observed that every
positive integer is the sum of four integer squares, but it was
not until 1770 that a complete proof of this remarkable fact was 
given by Lagrange. Also in 1770, Waring proposed a generalization 
of the four square theorem that became known as Waring's problem
and arguably led to the emergence of additive number theory. In 
modern terminology, Waring's conjecture states that for every 
integer $k \ge 2$ there exists an integer $s = s(k)$ such that 
every natural number $n$ is the sum of at most $s$ $k$th powers of 
natural numbers. Several special cases of this conjecture were 
settled during the 19th century, but the complete solution eluded 
mathematicians until 1909, when  Hilbert~\cite{Hilb09} proved the
existence of such an $s$ for all $k$ by means of a difficult 
combinatorial argument. 

Let $g(k)$ denote the least possible $s$ as above. Hilbert's 
method produced a very poor bound for $g(k)$. Using the circle 
method, Hardy and Littlewood were able to improve greatly on 
Hilbert's bound for $g(k)$. In fact, through the efforts of many
mathematicians, the circle method in conjunction with elementary 
and computational arguments has led to a nearly complete
evaluation of $g(k)$. In particular, we know that $g(k)$ is 
determined by certain special integers $n < 4^k$ that can only be 
represented as sums of a large number of $k$th powers of $1$, $2$ 
and $3$ (see \cite[\S 1.1]{Vaug97} for further details on $g(k)$).

A much more difficult question, and one that leads to a much deeper
understanding of the additive properties of $k$th powers, is that
of estimating the function $G(k)$, defined as the least $s$ such
that every \emph{sufficiently large} positive integer $n$ is the
sum of $s$ $k$th powers. This function was introduced by Hardy and
Littlewood~\cite{HaLi22}, who obtained the bound
\begin{equation}\label{1.8}
  G(k) \le (k - 2)2^{k - 1} + 5.
\end{equation}
In fact, they proved more than that. Let $I_{k, s}(n)$ denote
the number of solutions of the diophantine equation
\begin{equation}\label{1.9}
  x_1^k + x_2^k + \dots + x_s^k = n
\end{equation}
in $x_1, \dots, x_s \in \mathbb N$. Hardy and Littlewood showed 
that if $ s \ge (k-2) 2^{k-1} + 5 $, then                            
\begin{equation} \label{1.10}                                        
  I_{k, s}(n) \sim                                                      
    \frac {\Gamma^s \left( 1 + {\textstyle \frac 1k} \right)}{\Gamma\left(
\textstyle \frac sk \right)}
    \mathfrak S_{k, s}(n) n^{s/k - 1} \qquad \text{ as } \quad n \to \infty,        
\end{equation}            
where $\Gamma$ stands for Euler's gamma-function
and $\mathfrak S_{k, s}(n)$ is an absolutely convergent              
infinite series, called the \emph{singular series}, such that
\[
  \mathfrak S_{k, s}(n) \ge c_1(k, s) > 0.
\]

While the upper bound \eqref{1.8} represents a tremendous improvement
over Hilbert's result, it is still quite larger than the trivial
lower bound $G(k) \ge k + 1$.\footnote{ Let $X$ be large.
If $n \le X$, any solution of
\eqref{1.9} must satisfy $1 \le x_1, \dots, x_s \le X^{1/k}$. There
at most $X^{s/k}$ such $s$-tuples, which yield at most 
$(1/{s!} + o(1))X^{s/k}$
distinct sums $x_1^k + \dots + x_s^k$. Thus, when $s \le k$, there are 
not enough sums of $s$ $k$th powers to represent all the integers.}
During the mid-1930s I.~M.~Vinogradov introduced several refinements
of the circle method that allowed him to obtain a series of
improvements on \eqref{1.8} for large $k$. In their most elaborate
version, Vinogradov's methods yield a bound of the form\footnote{
In this and similar results appearing later, one can obtain an
explicit expressions in place of the $O$-terms, but those are too
complicated to state here.}
\[
  G(k) \le 2k(\log k + O(\log\log k)).
\]
First published by Vinogradov~\cite{IVin59} in 1959, this bound
withstood any significant improvement until 1992, when
Wooley~\cite{Wool92a} proved that
\[
  G(k) \le k(\log k + \log\log k + O(1)).
\]
The latter is the sharpest bound to date for $G(k)$ when $k$ is 
large. For smaller $k$, one can obtain better results
by using more specialized techniques (usually refinements of the
circle method). The best known bounds for $G(k)$, $3 \le k \le 20$,
are of the form $G(k) \le F(k)$, with $F(k)$ given by Table 1 below.

\medskip

\begin{center}
\small
\setlength{\tabcolsep}{3.5pt}
\renewcommand{\arraystretch}{1.5}
\begin{tabular}{cccccccccccccccccccc}
  \hline
  $k$    & 3 & 4 & 5  & 6  & 7  & 8  & 9  & 10 & 11 & 12 & 13 & 14 & 15  & 16  & 17 
& 18  & 19  & 20  \\
  \hline
  $F(k)$ & 7 & 16 & 17 & 24 & 33 & 42 & 50 & 59 & 67 & 76 & 84 & 92 & 100 & 109 &
117 & 125 & 134 & 142 \\
  \hline
\end{tabular}\\
$\phantom{.}$\\
{\sc Table 1.} Bounds for $G(k)$, $3 \le k \le 20$.
\end{center}

\medskip

\noindent
With the exception of the bound $G(3) \le 7$, all of these 
results have been obtained by an iterative version of the circle 
method that originated in the work of Davenport~\cite{Dave39a, 
Dave42} and Davenport and Erd\"os~\cite{DaEr39}. The bound for
$G(3)$ was established first by Linnik~\cite{Linn43} and until
recently lay beyond the reach of the circle method. The result 
on $G(4)$ is due to Davenport~\cite{Dave39b}, and in fact states 
that $G(4) = 16$. This is because $16$ biquadrates are needed to 
represent integers of the form $n = 31 \cdot 16^r$, $r \in 
\mathbb N$. Other than Lagrange's four squares theorem, this is 
the only instance in which the exact value of $G(k)$ is known.
However, Davenport~\cite{Dave39b} also proved that if $s \ge 14$,
all sufficiently large integers $n \equiv r \pmodulo {16}$, $1 
\le r \le s$, can be written as the sum of $s$ biquadrates; 
Kawada and Wooley \cite{KaWo99} obtained a similar result for as 
few as $11$ biquadrates. The remaining bounds in Table 1 appear 
in a series of recent papers by Vaughan and Wooley 
\cite{VaWo95}--\cite{VaWo00}.

A great deal of effort has also been dedicated to estimating the
function $\tilde G(k)$, which represents the least $s$ for which
the asymptotic formula \eqref{1.10} holds. For large $k$,
Ford~\cite{Ford95} showed that
\begin{equation}\label{1.11}
  \tilde G(k) \le k^2(\log k + \log\log k + O(1)),
\end{equation}
thus improving on earlier work by Vinogradov~\cite{IVin47},
Hua~\cite{Hua49}, and Wooley~\cite{Wool92b}. Furthermore,
Vaughan~\cite{Vaug86a, Vaug86b} and Boklan~\cite{Bokl94}
obtained the bounds
\[
  \tilde G(k) \le 2^k \quad (k \ge 3)
  \qquad \text{and} \qquad
  \tilde G(k) \le {\textstyle \frac 78} \cdot 2^k \quad (k \ge 6),
\]
which supersede \eqref{1.11} when $k \le 8$.

\medskip

The work on Waring's problem has inspired research on several other
questions concerned with the additive properties of $k$th powers 
(and of more general polynomial sequences). Such matters, however, 
are beyond the scope of this survey. The reader interested in a more
comprehensive introduction to Waring's problem should refer to the 
monographs \cite{ArChKa87, Vaug97} or to a recent survey article by
Vaughan and Wooley~\cite{VaWo02} (the latter also provides an 
excellent account of the history of Waring's problem).

\subsection{The Waring--Goldbach problem}
\label{s1.3}

Vinogradov's proof of the three prime theorem provided a blueprint
for subsequent applications of the Hardy--Littlewood circle method
to additive problems involving primes. Shortly after the publication 
of Theorem \ref{th1}, Vinogradov himself \cite{IVin37b} and
Hua~\cite{Hua38} began studying Waring's problem with prime variables,
known nowadays as the \emph{Waring--Goldbach problem}. They were able to 
generalize the asymptotic formula \eqref{1.3} to $k$th powers for
all $k \ge 1$ and ultimately their efforts led to the proof of 
Theorem \ref{th3} below.

In order to describe the current knowledge about the Waring--Goldbach
problem, we first need to introduce some notation. Let $k$ be a 
positive integer and $p$ a prime. We denote by $\theta = \theta 
(k, p)$ the (unique) integer such that $p^{\theta} \mid k$ and 
$p^{\theta+1} \nmid k$, and then define 
\begin{equation} \label{1.12}
  \gamma = \gamma (k, p) =
  \begin{cases}
    \theta + 2,  & \text{if } p=2, \; 2 \mid k, \\ 
    \theta + 1,  & \text{otherwise,}
  \end{cases} \qquad
  K(k) = \prod_{(p-1) \mid k} p^{\gamma}.
\end{equation}
In particular, we have $K(1) = 2$. It is not difficult to show that
if an integer $n$ is the sum of $s$ $k$th powers of primes greater
than $k + 1$, then $n$ must satisfy the congruence condition $n \equiv
s \pmodulo{K(k)}$. Furthermore, define
\begin{equation}\label{1.13}
  S(q, a) = \sum_{ \substack{ h = 1\\ (h, q) = 1}}^q
  e \left( \frac {ah^k}q \right), \quad 
  \mathfrak S_{k, s}^*(n) = \sum_{q = 1}^{\infty}
  \sum_{ \substack{ a = 1\\ (a, q) = 1}}^q
  \frac {S(q, a)^s}{\phi(q)^s} e\left( \frac {-an}q \right),
\end{equation}
where $(a, q)$ stands for
the greatest common divisor of $a$ and $q$, and $\phi(q)$ is Euler's
totient function, that is, the number of positive integers $n \le q$
which are relatively prime to $q$. The following result will be 
established in \S\ref{ssii2} and \S\ref{s3.4}.

\begin{theorem}\label{th3}
  Let $k, s$ and $n$ be positive integers, and let $R_{k, s}^*(n)$
  denote the number of solutions of the diophantine equation 
  \begin{equation}\label{1.14}
    p_1^k + p_2^k + \dots + p_s^k = n
  \end{equation}
  in primes $p_1, \dots, p_s$. Suppose that
  \[
    s \ge
    \begin{cases}
      2^k + 1,                          & \text{if } \;\; 1  \le k \le 5, \\
      \frac 78 \cdot 2^k + 1,           & \text{if } \;\; 6   \le k \le 8, \\
      k^2(\log k + \log\log k + O(1)),  \; & \text{if } \; \;  k > 8. 
    \end{cases}
  \]
  Then
  \begin{equation} \label{1.15}
    R_{k, s}^*(n) \sim \frac{\Gamma^s \left( 1 + {\textstyle \frac 1k} \right)}
    {\Gamma \left( \textstyle \frac sk \right)  }
    \mathfrak S_{k, s}^*(n) \frac {n^{s/k - 1}}{(\log n)^s} 
    \qquad \text{ as } \quad n \to \infty,
  \end{equation}
  where $\mathfrak S_{k, s}^*(n)$ is defined by \eqref{1.13}.
  Furthermore, the singular series $\mathfrak S_{k, s}^*(n)$ is
  absolutely convergent, and if $n \equiv s \pmodulo{K(k)}$, then
  $\mathfrak S_{k, s}^*(n) \ge c_2(k, s) > 0$.
\end{theorem}

In particular, we have the following corollaries to Theorem \ref{th3}.

\begin{corollary}\label{c31}
  Every sufficiently large integer $n \equiv 5 \pmodulo {24}$ can be
  represented as the sum of five squares of primes.
\end{corollary}

\begin{corollary}\label{c32}
  Every sufficiently large odd integer can be represented as the sum
  of nine cubes of primes.
\end{corollary}

Hua introduced a function $H(k)$ similar to the function $G(k)$ in 
Waring's problem. $H(k)$ is defined as the least integer $s$ such 
that equation \eqref{1.14} has a solution in primes $p_1, \dots, p_s$
for all sufficiently large $n \equiv s \pmodulo {K(k)}$. It is
conjectured that $H(k) = k + 1$ for all $k \ge 1$, but this
conjecture has not been proved for any value of $k$ yet. When
$k \le 3$, the sharpest known upper bounds for $H(k)$ are those 
given by Theorem~\ref{th3}, that is,
\[
  H(1) \le 3, \quad H(2) \le 5, \quad H(3) \le 9.
\]
When $k \ge 4$, the best results in the literature are as follows.

\begin{theorem}\label{th4}
  Let $k \ge 4$ be an integer, and let $H(k)$ be as above. Then
  \[
    H(k) \le
    \begin{cases}
      F(k),                 & \; \text{if } \;\; 4 \le k \le 10, \\
      k(4\log k + 2\log\log k + O(1)), & \; \text{if } \;\; k > 10, 
    \end{cases}
  \]
  where $F(k)$ is given by the following table.

  \begin{center}
    \small\rm
    \renewcommand{\arraystretch}{1.5}
    \begin{tabular}{cccccccc}
      \hline
      $k$    & 4  & 5  & 6  & 7  & 8  & 9  & 10  \\ \hline
      $F(k)$ & 14 & 21 & 33 & 46 & 63 & 83 & 107 \\ \hline
    \end{tabular}\\
    $\phantom{.}$\\
    {\sc Table 2.} Bounds for $H(k)$, $4 \le k \le 10$.
  \end{center}
\end{theorem}

The cases $k = 6$ and $8 \le k \le 10$ of Theorem \ref{th4} are 
due to Thanigasalam~\cite{Than87}, and the cases $k = 4, 5$ and 
$7$ are recent results of Kawada and Wooley~\cite{KaWo01} and 
Kumchev~\cite{Kumc03}, respectively. The bound for $k > 10$
is an old result of Hua, whose proof can be found in Hua's book 
\cite{Hua65}. To the best of our knowledge, this is the strongest
published result for large $k$, although it is well-known to 
experts in the field that better results are within the reach 
of Wooley's refinement of Vinogradov's methods. 
In particular, by inserting Theorem 1 in Wooley \cite{Wool95} into 
the machinery developed in Hua's monograph, one obtains
\[
  H(k) \le k({\textstyle\frac 32}\log k + O(\log\log k))
  \qquad \text {for } \;\; k \to \infty.
\]

\subsection{Other additive problems involving primes}                
\label{s1.4}

There are several variants and generalizations of the
Waring--Goldbach problem that have attracted a lot of attention
over the years. For example, one may consider the diophantine
equation
\begin{equation} \label{1.16}
  a_1p_1^k + a_2p_2^k + \dots + a_sp_s^k = n,
\end{equation}
where $n, a_1, \dots, a_s $ are fixed, not necessarily positive,
integers. There are several questions that we can ask about
equations of this form. The main question, of course, is that
of solubility. Furthermore, in cases where we do know that 
\eqref{1.16} is soluble, we may want to count the solutions with 
$p_1, \dots, p_s \le X$, where $X$ is a large parameter. A famous 
problem of
this type is the \emph{twin-prime conjecture}: there exist
infinitely many primes $p$ such that $p + 2$ is also prime, 
that is, the equation
\[
  p_1 - p_2 = 2
\]
has infinitely many solutions. It is believed that this conjecture
is of the same difficulty as the binary Goldbach problem, and in
fact, the two problems share a lot of common history. In particular,
while the twin-prime conjecture is still open, Chen's proof of
Theorem \ref{th2} can be easily modified to establish that there
exist infinitely many primes $p$ such that $p + 2 = P_2$.

Other variants of the Waring--Goldbach problem consider more
general diophantine equations of the form
\[
  f(p_1) + f(p_2) + \dots + f(p_s) = n,
\]
where $f(X) \in \mathbb Z[X]$, or systems of equations of the
types \eqref{1.1} or \eqref{1.16}. For example, Chapters 10 and 11
in Hua's monograph~\cite{Hua65} deal with the system
\[
  p_1^j + p_2^j + \dots + p_s^j = n_j \qquad (1 \le j \le k).
\]
The number of solutions of this system satisfies an asymptotic
formula similar to \eqref{1.15}, but the main term in that
asymptotic formula is less understood than the main term in
\eqref{1.15} (see \cite{ArCh02a, Ch86a, Hua65, Mitk92} for 
further details).

Another classical problem in which a system of diophantine equations
arises naturally concerns the existence of non-trivial arithmetic 
progressions consisting of $r$ primes. It has been conjectured that
for every integer $r \ge 3$ there are infinitely many such arithmetic 
progressions. In other words, the linear system
\[
 p_i - 2 p_{i+1} + p_{i+2} = 0  \qquad (1 \le i \le r - 2)
\]
has infinitely many solutions in distinct primes $p_1, \dots, p_r$.
In the case $r = 3$ this can be established by a variant of 
Vinogardov's proof of the three primes theorem, but when $r > 3$ the 
above system lies beyond the reach of the circle method. In fact, 
until recently the most significant insight into progressions 
of more than three primes were the following two results:
\begin{itemize}
  \item Heath-Brown \cite{H-Br81} succeeded to prove that there 
  exist infinitely many arithmetic progressions of three primes and 
  a $P_2$-number.
  \item Balog~\cite{Balog92} proved that for any $r$ there are $r$ 
  distinct primes $p_1, \dots, p_r$ such that all the averages 
  $\frac{1}{2}(p_i + p_j)$ are prime.
\end{itemize}
Thus, the specialists in the field were stunned when Green and Tao 
\cite{GrTa04} announced their amazing proof of the full conjecture.
The reader will find a brief description of their ideas and of some
related recent work in the last section.

Finally, instead of \eqref{1.1}, one may study the inequality
\[
  |x_1 + \dots + x_s - \alpha| < \eps,
\]
where $\alpha$ is a real number,
$\eps$ is a small positive number and
$x_1, \dots, x_s$ are real variables taking values from a
given sequence (or sequences). For example, by setting
$x_j = p_j^c$ where $c > 1$ in not an integer, we can
generalize the Waring--Goldbach problem to fractional powers
of primes. We will mention several results of this form in
\S \ref{DI}.

\section{The distribution of primes}
\label{s2}

In this section we discuss briefly some classical results about
primes, which play an important role in additive prime number 
theory.

\subsection{The Prime Number Theorem}
\label{s2.1}

The first result on the distribution of primes is Euclid's
theorem that there are infinitely many prime numbers. In 1798
Legendre conjectured that the prime counting function $\pi(x)$
(i.e., the number of primes $p \le x$) satisfies the asymptotic      
relation
\begin{equation}\label{2.1}
  \lim_{x \to \infty} \frac {\pi(x)}{x/(\log x)} = 1;
\end{equation}
this is the classical statement of the Prime Number Theorem.         
Later Gauss observed that the logarithmic integral
\[
  \li x = \int_2^x \frac{dt}{\log t}
\]
seemed to provide a better approximation to $\pi(x)$ than the
function $x/(\log x)$ appearing in \eqref{2.1}, and this is 
indeed the case. Thus, in anticipation of versions of the Prime 
Number Theorem that are more precise than \eqref{2.1}, we define 
the error term
\begin{equation}\label{2.2}                                          
  \Delta(x) = \pi(x) - \li x.                                        
\end{equation}                                                       

The first step toward a proof of the Prime Number Theorem was 
made by Chebyshev. In the early 1850s he proved that 
\eqref{2.1} predicts correctly the order of $\pi(x)$, that is,
he established the existence of absolute constants $c_2 > c_1
> 0$ such that
\[
  \frac{c_1x}{\log x} \le \pi(x) \le \frac{c_2x}{\log x}.
\]
Chebyshev also showed that if the limit on the left side of
\eqref{2.1} exists, then it must be equal to $1$.

In 1859 Riemann published his famous memoir \cite{Riem59}, in
which he demonstrated the intimate relation between $\pi(x)$
and the function which now bears his name, that is, the
{\em Riemann zeta-function} defined by
\begin{equation}\label{2.3}
  \zeta(s) = \sum_{n = 1}^{\infty} n^{-s}
  = \prod_p \left( 1 - p^{-s} \right)^{-1}
  \quad (\re(s) > 1).
\end{equation}
This and similar series had been used earlier by Euler\footnote{     
In particular, Euler established the equality between $\zeta(s)$     
and the infinite product in \eqref{2.3}, which is known as the       
\emph{Euler product} of $\zeta(s)$.} and Dirichlet, but only as      
functions of a real variable. Riemann observed that $\zeta(s)$
is holomorphic in the half-plane $\re(s) > 1$ and that it can be
continued analytically to a meromorphic function, whose only
singularity is a simple pole at $s = 1$. It is not difficult to
deduce from \eqref{2.3} that $\zeta(s) \ne 0$ in the half-plane
$\re(s) > 1$. Riemann observed that $\zeta(s)$ has infinitely
many zeros in the strip $0 \le \re(s) \le 1$ and proposed several 
conjectures concerning those zeros and the relation between them 
and the Prime Number Theorem. The most famous among those 
conjecture---and the only one that is still open---is known as 
the {\em Riemann Hypothesis}. 

\begin{RH}[RH]
  All the zeros of $\zeta(s)$ with $0 \le \re(s) \le 1$ lie on 
  the line $\re(s) = \frac 12$.
\end{RH} 

The remaining conjectures in Riemann's paper were proved by the
end of the 19th century. In particular, it was proved that the
Prime Number Theorem follows from the nonvanishing of $\zeta(s)$ 
on the line $\re(s) = 1$. Thus, when in 1896 Hadamard and de la 
Vall\'ee Poussin proved (independently) 
that $\zeta(1 + it) \ne 0$ for all real $t$, the Prime Number 
Theorem was finally proved. In 1899 de la Vall\'ee Poussin  
obtained the following quantitative result.\footnote{
Functions 
of the type $f(x) = \exp \big( (\log x)^{\lambda} \big)$, where 
$\lambda$ is a constant, are quite common in analytic number
theory. To help the reader appreciate results such as Theorems 
\ref{PNT} and \ref{th6}, we remark that as $x \to \infty$ such
a function with $0 < \lambda < 1$ grows more rapidly than any
fixed power of $\log x$, but less rapidly than $x^{\eps}$ 
for any fixed $\eps > 0$.} (Henceforth, we 
often use Vinogradov's notation $A \ll B$, which means that 
$A = O(B)$.)                    

\begin{theorem}[de la Vall\'ee Poussin, 1899]\label{PNT}
  Let $\Delta(x)$ be defined by \eqref{2.2}. There exists an
  absolute constant $c > 0$ such that
  \[
    \Delta(x) \ll x \exp \big( -c \sqrt{\log x} \big).
  \]
\end{theorem}

De la Vall\'ee Poussin's theorem has been improved somewhat, but
not nearly as much as one would hope. The best result to date is
due to I. M. Vinogradov~\cite{IVin58} and Korobov~\cite{Koro58},
who obtained (independently) the following estimate for $\Delta(x)$.

\begin{theorem}[Vinogradov, Korobov, 1958]\label{th6}
  Let $\Delta(x)$ be defined by \eqref{2.2}. There exists an
  absolute constant $c > 0$ such that
  \[
    \Delta(x) \ll x \exp \big( -c (\log x)^{3/5}
  (\log \log x)^{-1/5} \big).
  \]
\end{theorem}

In comparison, if the Riemann Hypothesis is assumed, one has
\begin{equation}\label{2.4}
  \Delta(x) \ll x^{1/2}\log x,
\end{equation}
which, apart from the power of the logarithm, is best possible.
The reader can find further information about the Prime Number
Theorem and the Riemann zeta-function in the standard texts on
the subject (e.g., \cite{Dave00, Ivic85, KaVo92, Prac57, Titc86}).

\subsection{Primes in arithmetic progressions}
\label{ssiii2}

In a couple of memoirs published in 1837 and 1840, Dirichlet 
proved that if $a$ and $q$ are natural numbers with $(a, q) = 1$,
then the arithmetic progression $a \; \modulo q$ contains 
infinitely many primes. In fact, Dirichlet's argument
can be refined as to establish the asymptotic formula
\begin{equation}\label{2.5}
  \sum_{\substack{ p \le x \\ p \equiv a \pmodulo q}}
  \frac {\log p}{p} \sim \frac {1}{\phi(q)} 
  \sum_{p \le x} \frac {\log p}{p}
  \qquad \text{as } \quad x \to \infty,
\end{equation}
valid for all $a$ and $q$ with $(a, q) = 1$. Fix $q$ and consider
the various arithmetic progressions $a \; \modulo q$ (here $\phi(q)$
is Euler's totient function). Since all but finitely many primes 
lie in progressions with $(a, q) = 1$ and there are $\phi(q)$ such 
progressions, \eqref{2.5} suggests that each arithmetic
progression $a \; \modulo q$, with $(a, q) = 1$, ``captures its
fair share'' of prime numbers, i.e., that the primes are uniformly
distributed among the (appropriate) arithmetic progressions to a
given modulus $q$. Thus, one may expect that if $(a, q)=1$, then     
\begin{equation}\label{2.6}
  \pi(x; q, a) =
  \sum_{\substack{ p \le x \\ p \equiv a \pmodulo q}} 1
  \sim \frac {\li x}{\phi(q)} \qquad \text{as } \quad  x \to \infty.
\end{equation}
This is the \emph{prime number theorem for arithmetic progressions}. 
One may consider \eqref{2.6} from two different view points. First,
one may fix $a$ and $q$ and ask whether \eqref{2.6} holds (allowing
the convergence to depend on $q$ and $a$). Posed in this form, the
problem is a minor generalization of the Prime Number Theorem. In
fact, shortly after proving Theorem \ref{PNT}, de la Vall\'ee
Poussin established that
\[
  \Delta(x; q, a) = \pi(x; q, a) - \frac {\li x}{\phi(q)} \ll
  x \exp \big( -c \sqrt{\log x} \big),
\]
where $c = c(q, a) > 0$ and the implied constant depends on $q$
and $a$. The problem becomes much more difficult if one requires
an estimate that is explicit in $q$ and uniform in $a$. The first
result of this kind was obtained by Page~\cite{Page35}, who 
proved the existence of a (small) positive number $\delta$ such 
that
\begin{equation}\label{2.7}
  \Delta(x; q, a) \ll x \exp \left( -(\log x)^{\delta} \right),
\end{equation}
whenever $1 \le q \le (\log x)^{2 - \delta}$ and $(a, q) = 1$.
In 1935 Siegel~\cite{Sieg35} (essentially) proved the following
result known as the Siegel--Walfisz theorem.

\begin{theorem}[Siegel, 1935]\label{SW}
  For any fixed $A > 0$, there exists a constant $c = c(A) > 0$ 
  such that
  \[ 
    \pi(x; q, a)= \frac {\li x}{\phi(q)} 
    + O\big( x \exp\big( -c \sqrt{\log x} \big) \big)
  \]
  whenever $q \le (\log x)^A$ and $(a,q)=1$.
\end{theorem}

\begin{remark}
  While this result is clearly sharper than Page's, it does have
  one significant drawback: it is ineffective, that is, given
  a particular value of $A$, the proof does not allow the
  constant $c(A)$ or the $O$-implied constant to be computed.
\end{remark}

The above results have been proved using the analytic properties
of a class of generalizations of the Riemann zeta-function known as
\emph{Dirichlet $L$-functions}. For each positive integer $q$ there
are $\phi(q)$ functions $ \chi : \mathbb Z \to \mathbb C $, called
{\em Dirichlet characters $\modulo q$}, with the following           
properties:                                                          
\begin{itemize}                                                      
  \item $\chi$ is {\em totally multiplicative}:                      
    $\chi(mn) = \chi(m)\chi(n)$;                                    
  \item $\chi$ is $q$-periodic;                        
  \item $|\chi(n)| = 1$ if $(n, q) = 1$ and $\chi(n) = 0$ if       
    $(n, q) > 1$;
  \item if $(n, q) = 1$, then
    \[
      \sum_{\chi \; \modulo q} \chi(n) 
      = \begin{cases}
          \phi(q)  & \text{ if } \;\; n \equiv 1 \pmodulo {q}, \\
          0        & \text{ otherwise}.
      \end{cases}                                            
    \]
\end{itemize}                 
For more information about the construction and properties of the 
Dirichlet characters we refer the reader to \cite{Dave00, IwKo04, 
Kara93, Prac57}.

Given a character $\chi$ $\modulo q$, we define the
Dirichlet $L$-function
\[
  L(s, \chi) = \sum_{n = 1}^{\infty} \chi(n)n^{-s} =
  \prod_p \left( 1 - \chi(p)p^{-s} \right)^{-1}
  \quad (\re(s) > 1).
\]
Similarly to $\zeta(s)$, $L(s, \chi)$ is holomorphic in the          
half-plane $\re(s) > 1$ and can be continued analytically to a       
meromorphic function on $\mathbb C$ that has at most one pole,
which (if present) must be a simple pole at $s = 1$. Furthermore,
just as $\zeta(s)$, the continued $L(s, \chi)$ has infinitely many
zeros in the strip $0 \le \re(s) \le 1$, and the horizontal
distribution of those zeros has important implications on the
distribution of primes in arithmetic progressions. For example,
the results of de la Vall\'ee Poussin, Page and Siegel mentioned
above were proved by showing that no $L$-function can have a zero
``close" to the line $\re(s) = 1$. We also have the following
generalization of the Riemann Hypothesis.

\begin{GRH}[GRH]
  Let $L(s, \chi)$ be a Dirichlet $L$-function.
  Then all the zeros
  of $L(s, \chi)$ with $0 \le \re(s) \le 1$ lie on the line 
  $\re(s) = \frac 12$.
\end{GRH} 

Assuming GRH, we can deduce easily that if $(a, q) = 1$, then
\begin{equation}\label{2.8}
  \pi(x; q, a)= \frac {\li x}{\phi(q)}
  + O\left( \textstyle x^{1/2}\log x \right),
\end{equation}
which is nontrivial when $1 \le q \le x^{1/2}(\log x)^{-2 - \eps}$.

In many applications one only needs \eqref{2.8} to hold ``on average'' 
over the moduli $q$. During the 1950s and 1960s several authors 
obtained estimates for averages of $\Delta(x; q, a)$. In particular, 
the following quantity was studied extensively:
\[
  E(x, Q) = \sum_{q \le Q} \max_{(a, q) = 1} \max_{y \le x}
  |\Delta(y; q, a)|.
\]
The trivial bound for this quantity is $E(x, Q) \ll x \log x$.
One usually focuses on finding the largest value of $Q$ for which 
one can improve on this trivial bound, even if the improvement is 
fairly modest. The sharpest result in this direction was 
established (independently) by Bombieri~\cite{Bomb65} and
A.~I.~Vinogradov~\cite{AVin65} in 1965. Their result is known as
the Bombieri--Vinogradov theorem and (in the slightly stronger
form given by Bombieri) can be stated as follows.

\begin{theorem}[Bombieri, Vinogradov, 1965]\label{B-V}
  For any fixed $A > 0$, there exists a $B = B(A) > 0$ such that
  \begin{equation}\label{iii7}
    E(x, Q) \ll x (\log x)^{-A},
  \end{equation}
  provided that $Q \le x^{1/2}(\log x)^{-B}$.
\end{theorem}

We should note that other than the value of $B(A)$ the range for
$Q$ in this result is as long as the range we can deduce from GRH.
Indeed, GRH yields $B = A + 1$, whereas Bombieri obtained Theorem
\ref{B-V} with $B = 3A + 22$ and more recently Vaughan~\cite{Vaug75}
gave $B = A + 5/2$.

\subsection{Primes in short intervals}
\label{ssiii3}

Throughout this section, we write $p_n$ for the $n$th prime number.
We are interested in estimates for the difference $p_{n + 1} - p_n$
between two consecutive primes.
Cram\'er was the first to study this question systematically. He 
proved \cite{Cram20} that the Riemann Hypothesis implies 
\[
  p_{n + 1} - p_n \ll p_n^{1/2}\log p_n.
\]
Cram\'er also proposed a probabilistic model of the prime numbers
that leads to very precise (and very bold) predictions of the 
asymptotic properties of the primes. In particular, he conjectured 
\cite{Cram37} that 
\begin{equation}\label{iii10}
  \limsup_{n \to \infty} \frac{p_{n + 1} - p_n}{\log^2 p_n} = 1.
\end{equation}

A non-trivial upper bound for $p_{n + 1} - p_n$                           
can be obtained as a consequence of the                              
Prime Number Theorem, but Hoheisel~\cite{Hohe30}                    
found a much sharper result. He proved unconditionaly                 
the asymptotic formula                                              
\begin{equation}\label{iii8}
  \pi(x + h) - \pi(x) \sim h(\log x)^{-1} \qquad \text{as } \quad
  x \to \infty,
\end{equation}
with $h = x^{1 - (3300)^{-1}}$. There have been several improvements 
on Hoheisel's result and it is now known that \eqref{iii8} holds
with $h = x^{7/12}$ (see Heath-Brown~\cite{H-Br88}).
Furthermore, several mathematicians have shown that
even shorter intervals must contain primes (without establishing
an asymptotic formula for the number of primes in such intervals).
The best result in this directions is due to
Baker, Harman, and Pintz~\cite{BaHaPi01}, who proved that
for each $n$ one has
\[
  p_{n + 1} - p_n \ll p_n^{0.525}.
\]

A related problem seeks small gaps between consecutive primes.
In particular, the twin-prime conjecture can be stated in the form
\[
  \liminf_{n \to \infty} (p_{n + 1} - p_n) = 2.
\]
It is an exercise to show that the Prime Number Theorem
implies the inequality
\[
  \liminf_{n \to \infty} \frac {p_{n + 1} - p_n}{\log p_n} \le 1.
\]
Improvements on this trivial bound, on the other hand, have proved
notoriously difficult and, so far, the best result, due to Maier
\cite{Maie88}, is
\[
  \liminf_{n \to \infty} \frac {p_{n + 1} - p_n}{\log p_n} \le 0.2486\dots.
\]

\subsection{Primes in sparse sequences}
\label{ssiii4}

We say that an infinite sequence of primes $\mathcal S$ is
{\it sparse} if 
\[
  \pi(\mathcal S; x) := \# \big\{ p \in \mathcal S : p \le x \big \}
  = o(\pi(x)) \qquad \text{as } \quad x \to \infty.
\]
A classical example that has attracted a great deal of attention
but has proved notoriously difficult is that of primes represented
by polynomials. To this day, there is not a single example of a 
polynomial $f(X) \in \mathbb Z[X]$ of degree at least $2$ which is 
known to take on infinitely many prime values. The closest 
approximation is a result of Iwaniec \cite{Iwan78}, who showed 
that if $a, b, c$ are integers such that $a > 0$, $(c, 2) = 1$, 
and the polynomial $f(X) = aX^2 + bX + c$ is irreducible, then 
$f(X)$ takes on infinitely many $P_2$-numbers. On the other hand,
in recent years there has been some exciting progress in the 
direction of finding polynomials in two variables that represent 
infinitely many primes. In 1998 Friedlander and 
Iwaniec~\cite{FrIw98a} proved that the polynomial $X^2 + Y^4$ 
represents infinitely many primes. We note that this polynomial 
takes on $O(x^{3/4})$ values up to $x$. In 2001 
Heath-Brown~\cite{H-Br01} obtained an analogous result for the 
polynomial $X^3 + 2Y^3$ whose values are even sparser: it takes 
on $O(x^{2/3})$ values up to $x$. Furthermore, Heath-Brown and 
Moroz~\cite{H-BMo02} extended the latter result to general 
irreducible binary cubic forms in $\mathbb Z[X, Y]$ (subject to 
some mild necessary conditions).

Another class of sparse sequences of prime numbers arises in the
context of diophantine approximation. The two best known examples
of this kind are the sequences
\begin{equation} \label{thinset1}
  \mathcal S_{\lambda} = \left\{ p : p \text{ is prime with }
  \{ \sqrt{p} \} < p^{ - \lambda } \right\}
\end{equation}
and
\begin{equation}\label{thinset2}
  \mathcal P_c = \left\{ p : p = [n^c]
  \text{ for some integer } n \right \}.
\end{equation}
Here, $\lambda \in (0, 1)$ and $c > 1$ are fixed real numbers, 
$\{ x \}$ denotes the fractional part of the real number $x$, 
and $[x] = x - \{ x \}$. The sequence $\mathcal S_{\lambda}$ was
introduced by I.~M.~Vinogradov, who proved (see 
\cite[Chapter~4]{IVin76}) that if $0 < \lambda < 1/10$, then
\[
    \pi(\mathcal S_{\lambda}; x) \sim
    \frac{ x^{ 1 - \lambda } } { (1- \lambda) \log x} \qquad
     \text{ as } \quad x \to \infty.
\]
The admissible range for $\lambda$ has
been subsequently extended to $0 < \lambda < 1/4$ by Balog 
\cite{Balog83} and Harman \cite{Harm83b}, while Harman and Lewis
\cite{HaLe01} showed that $\mathcal S_{\lambda}$ is infinite for
$0 < \lambda < 0.262$. 

The first to study the sequence \eqref{thinset2} was 
Piatetski-Shapiro \cite{PS53}, who considered $\mathcal P_c$ as
a sequence of primes represented by a ``polynomial of degree
$c$''. Piatetski-Shapiro proved that $\mathcal P_c$
is infinite when $1 < c < 12/11 $. The range for $c$ has
been extended several times and it is currently known (see
Rivat and Wu~\cite{RiWu01}) that $\mathcal P_c$ is infinite
when $1 < c < 243/205$. Furthermore, it is known (see Rivat and
Sargos \cite{RiSa01}) that when $1 < c < 1.16117\dots$, we have
\[
  \pi(\mathcal P_c; x) \sim \frac {x^{1/c}}{\log x} 
      \qquad \text{ as } \quad x \to \infty.
\]

\section{The Hardy--Littlewood circle method}
\label{sii}

Most of the results mentioned in the Introduction have been proved
by means of the Hardy--Littlewood circle method. In this section,
we describe the general philosophy of the circle method, using its
applications to the Goldbach and Waring--Goldbach problems to
illustrate the main points.

\subsection{Vinogradov's three prime theorem}
\label{ssii1}

\subsubsection{Preliminaries}
\label{dl1} 

Using the orthogonality relation
\begin{equation} \label{ii1}
  \int_0^1 e(\alpha m) \, d \alpha =
  \begin{cases}
    1, &  \text{if } \; m = 0, \\
    0, &  \text{if } \; m \ne 0,
  \end{cases}
\end{equation}
we can express $R(n)$ as a Fourier integral. We have
\begin{align}\label{ii2}
  R(n) &= \sum_{p_1, p_2, p_3 \le n} \int_0^1
  e \left( \alpha \left( p_1 + p_2 + p_3 - n \right) \right) \, d\alpha \\
  &= \int_0^1 f(\alpha)^3 e(-\alpha n) \, d\alpha, \notag
\end{align}
where $f(\alpha)$ is the exponential sum \eqref{1.2}.

The circle method uses \eqref{ii2} to derive an asymptotic formula 
for $R(n)$ from estimates
for $f(\alpha)$. The analysis of the right side of \eqref{ii2} rests
on the observation that the behavior of $f(\alpha)$ depends on the 
distance from $\alpha$ to the set of fractions with ``small'' 
denominators. When $\alpha$ is ``near'' such a fraction, we expect
$f(\alpha)$ to be ``large'' and to have certain asymptotic behavior.
Otherwise, we can argue that the numbers $e(\alpha p)$ are 
uniformly distributed on the unit circle and hence $f(\alpha)$ 
is ``small". In order to make these observations rigorous, we need  
to introduce some notation. Let $B$ be a positive constant
to be chosen later and set
\begin{equation} \label{ddd1}
   P = (\log n)^B.
\end{equation}
If $a$ and $q$ are integers with $1 \le a \le q \le P $                    %7dec
and $(a, q)=1$, we define the \emph{major arc}\footnote{This term 
may seem a little 
peculiar at first, given that $\mathfrak M(q, a)$ is in fact an 
interval. The explanation is that, in the original version of the 
circle method, Hardy and Littlewood used power series and Cauchy's
integral formula instead of exponential sums and \eqref{ii1} (see 
\cite[\S 1.2]{Vaug97}). In that setting, the role of $\mathfrak 
M(q, a)$ is played by a small circular arc near the root of unity 
$e(a/q)$.}
\begin{equation}\label{ii3}
  \mathfrak M(q, a) =                                                      %7dec
       \left[ \frac{a}{q} - \frac{P}{qn} \, ,                              %7dec
             \, \frac{a}{q} + \frac{P}{qn}  \right] .                      %7dec 
%\left\{ \alpha \in [0, 1] :                                               %7dec
%  |q\alpha - a| < Pn^{-1} \right\}.                                       %7dec 
\end{equation}                                                             %7dec 
The integration in \eqref{ii2} can be taken over any interval              %7dec  
of length one and, in particular, over                                     %7dec  
$ \big[ Pn^{-1} , 1 + Pn^{-1} \big] $.                                     %7dec 
We partition this interval into two subsets:                               %7dec 
\begin{equation}\label{ii4}                                                %7dec 
  \mathfrak M = \bigcup_{q \le P} \bigcup_{ \substack{1 \le a \le q\\      %7dec 
  (a, q) = 1}} \mathfrak M(q, a) \quad \text{and} \quad                    %7dec 
  \mathfrak m = \big[ Pn^{-1}, 1 + Pn^{-1} \big] \setminus \mathfrak M,    %7dec 
\end{equation}                                                             %7dec
%                                                                          %7dec
%   TUK NAPRAVIH MALKA PROMJANA, PONEZHE INACHE IMA                        %7dec
%   LEKO RAZMINAWANE MEZHDU FORMULI (3.4) I (3.7) - (3.8).                 %7dec
%   STRUWA MI SE CHE TOVA, KOETO SYM NAPRAVIL                              %7dec
%   E NAJ-BEZBOLEZNENO, A I SYSHTO E SHIROKO RAZPROSTRANENO.               %7dec
%   NO AKO NE TI HARESVA - PROMENI GO.                                     %7dec
%                                                                          %7dec
called respectively the \emph{set of major arcs} and the \emph{set of
minor arcs}. Then from \eqref{ii2} and \eqref{ii4} it follows that
\begin{equation}\label{ii5}
  R(n) = R(n, \mathfrak M) + R(n, \mathfrak m),
\end{equation}
where we have denoted
\[
  R(n, \mathfrak B) = \int_{\mathfrak B} 
  f(\alpha)^3 e(-\alpha n) \, d\alpha.
\]
In the next section we explain how, using Theorem~\ref{SW} and 
standard results from elementary number theory,
one can obtain an asymptotic formula for $R(n, \mathfrak M)$ 
(see \eqref{ii15} below). Then in \S \ref{minorarcs1} and \S 
\ref{ssii1.4} we discuss how one can show that $R(n, \mathfrak m)$ 
is of a smaller order of magnitude than the main term in that 
asymptotic formula (see \eqref{ii16}).

\subsubsection{The major arcs}
\label{ssii1.2}

In this section we sketch the estimation of the contribution from
the major arcs. The interested reader will find the missing details
in \cite[Chapter 10]{Kara93} or \cite[Chapter 2]{Vaug97}. 

It is easy to see that the major arcs $\mathfrak M(q, a)$ are mutually 
disjoint. Thus, using \eqref{ii3} and \eqref{ii4}, we can write
\begin{equation} \label{5oct1}
  R(n, \mathfrak M) =   \sum_{q \le P} \,
     \sum_{ \substack{1 \le a \le q \\ (a, q) = 1 }} \,
  \int_{-P/(qn)}^{P/(qn)} 
    f ( a/q + \beta )^3 e \big(  -(a/q + \beta) n \big) \, d\beta.
\end{equation}

We now proceed to approximate $ f \big( a/q + \beta \big) $ by a simpler expression.
To motivate our choice of the approximation, we first consider the case $\beta = 0 $.
We split the sum $ f \big( a/q \big) $ into
subsums according to the residue of $p$ modulo $q$ and take into account
the definition \eqref{2.6}. We get                                                  
            
\[
  f \left( \frac aq \right)
  = \sum_{h = 1}^q \sum_{\substack{p \le n\\ p \equiv h \pmodulo{q}}} 
  e \left( \frac {ap}{q} \right) 
  = \sum_{h = 1}^q e \left( \frac {ah}{q} \right) \pi(n; q, h).
\]
The contribution of the terms with $(h, q)>1$ is negligible (at most $q$). If $(h,
q) = 1$, our choice \eqref{ddd1} of the parameter $P$ ensures that we can appeal to
Theorem \ref{SW} to approximate $\pi (n; q, h)$ by $\phi(q)^{-1}\li n$. We deduce
that 
\begin{equation}\label{ii7}
  f \left( \frac aq \right)
  = \frac {\li n}{\phi(q)} 
  \sum_{ \substack{ h = 1\\ (h, q) = 1}}^q 
  e \left( \frac {ah}{q} \right) + O \left( qnP^{-4} \right).
\end{equation}
The exponential sum on the right side of \eqref{ii7} is known as the
{\em Ramanujan sum} and is usually denoted by $c_q(a)$. Its value is
known for every pair of integers $a$ and $q$ (see \cite[Theorem 271]
{HaWr79}). In particular, when $(a, q) = 1$ we have
$c_q(a) = \mu(q)$, where $\mu$ is the M\"obius function
\begin{equation}\label{Moebius}
  \mu(n) = \begin{cases}
    1,      & \text{if }  n = 1, \\
    (-1)^k, & \text{if }  n = p_1 \cdots p_k
              \text{ is the product of } k \text{ distinct primes}, \\
    0,      & \text{otherwise}.
  \end{cases}
\end{equation}
The situation does not change much if instead of $\alpha = a/q$ we 
consider $\alpha = a/q + \beta \in \mathfrak M(q, a)$. In this case
we find that
\begin{equation}\label{ii9}
  f \left( \frac aq + \beta \right) = 
  \frac {\mu(q)}{\phi(q)} \cdot v(\beta) 
  + O \left( nP^{-3} \right),
\end{equation}
where
\[
  v(\beta) = \int_2^n \frac {e(\beta u)}{\log u} \, du.
\]

Raising \eqref{ii9} to the third power and inserting the result into 
the right side of \eqref{5oct1}, we obtain
\begin{equation}\label{ii10}
  R(n, \mathfrak M) = 
  \sum_{q \le P} \frac {\mu(q)c_q(-n)}{\phi(q)^3} 
  \int_{-P/(qn)}^{P/(qn)} v(\beta)^3 e( -\beta n) \, d\beta 
  + O \left( n^2P^{-1} \right). 
\end{equation}
At this point, we extend the integration over $\beta$ to the whole real line, and
then the summation over $q$ to all positive integers. The arising error terms can be
controlled easily by means of well-known bounds for the functions $v(\beta)$ and
$\phi(q)$, and we find that
\begin{equation}\label{ii12}
  R(n, \mathfrak M) = \mathfrak S(n) J(n) 
  + O \left( n^2P^{-1} \right),
\end{equation}
where $\mathfrak S(n)$ and $J(n)$ are the \emph{singular
series} and the \emph{singular integral} defined by
\begin{gather*}
  \mathfrak S(n) =
  \sum_{q = 1}^{\infty} \frac {\mu(q)c_q(-n)}{\phi(q)^3}, \qquad
  J(n) = \int_{-\infty}^{\infty} v(\beta)^3 e(-\beta n) \, d\beta.
\end{gather*}

The series $\mathfrak S(n)$ actually satisfies \eqref{1.4}.
Indeed, the function 
\[
  g(q) = \mu(q)c_q(-n)\phi(q)^{-3}
\] 
is {\em multiplicative} in $q$, that is, 
$g(q_1 q_2) = g(q_1)g(q_2)$ whenever $(q_1, q_2) = 1$. Hence, 
using the absolute convergence of $\mathfrak S(n)$ and the
elementary properties of the arithmetic functions involved
in the definition of $g(q)$, we can represent the singular
series as an Euler product:
\begin{align}
  \mathfrak S(n) =
  \sum_{q = 1}^{\infty} g(q)
  &= \prod_p \left( 1 + g(p) + g(p^2) + \cdots \right) 
          \notag \\
  &= \prod_{p \mid n} \left( 1 - \frac {1}{(p-1)^2} \right)
    \prod_{p \nmid n} \left( 1 + \frac {1}{(p-1)^3} \right).
    \notag
\end{align}
Also, an application of Fourier's inversion formula and some calculus
reveal that 
\[
  J(n) = \frac {n^2}{2(\log n)^3} 
  + O \left( n^2(\log n)^{-4} \right).
\]
Therefore, if $B \ge 4$ we can conclude that
\begin{equation}\label{ii15}
  R(n, \mathfrak M) = \frac {n^2}{2(\log n)^3} \mathfrak S(n)
  + O \left( n^2(\log n)^{-4} \right).
\end{equation}

\subsubsection{The minor arcs}
\label{minorarcs1}

In view of \eqref{ii5} and \eqref{ii15}, it suffices to prove that
(for some $B \ge 4$)
\begin{equation}\label{ii16}
  R(n, \mathfrak m) \ll n^2 (\log n)^{-4}.
\end{equation}
We have
\begin{equation}\label{ii17}
  |R(n, \mathfrak m)| \le \int_{\mathfrak m} |f(\alpha)|^3 \, d\alpha
  \le \Big( \sup_{\mathfrak m} |f(\alpha )| \Big) 
  \int_0^1 |f(\alpha)|^2 \, d\alpha.  
\end{equation}
By Parseval's identity and the Prime Number Theorem,
\begin{equation}\label{ii18}
  \int_0^1 |f(\alpha)|^2 \, d\alpha 
  = \sum_{p \le n} 1 \ll n(\log n)^{-1}.
\end{equation}
Thus, \eqref{ii16} will follow from \eqref{ii17}, if we show that
\begin{equation}\label{ii19}
  \sup_{\mathfrak m} |f(\alpha)| \ll n(\log n)^{-3}.
\end{equation}

We note that the trivial estimate for $f(\alpha)$ is
\[
   f(\alpha) \ll \sum_{ p \le n } 1 \ll n (\log n)^{-1} ,
\]
so in order to establish \eqref{ii19}, we have 
to save a power of $\log n$ over this trivial estimate
(uniformly with respect to $\alpha \in \mathfrak m $).
We can do this using the following lemma, which provides
such a saving under the assumption that $\alpha$ can be approximated 
by a reduced fraction whose denominator $q$ is ``neither too small,
nor too large.''

\begin{lemma}\label{lii2}
  Let $\alpha$ be real and let $a$ and $q$ be integers satisfying
  \[
    1 \le q \le n, \quad (a, q) = 1, \quad
    |q\alpha - a| \le q^{-1}.
  \]
  Then
  \[
    f(\alpha) \ll (\log n)^3 
    \left( nq^{-1/2} + n^{4/5} + n^{1/2}q^{1/2} \right).
  \]
\end{lemma}

This is the sharpest known version of the estimate for 
$f(\alpha)$ established by I. M. Vinogradov \cite{IVin37a} in 1937.
As we mentioned in the Introduction, that result was the main 
innovation in Vinogradov's proof of Theorem~\ref{th1}. The above
version is due to Vaughan~\cite{Vaug77}. 

We shall explain the proof of 
Lemma~\ref{lii2} in the next section and now we shall use it
to establish \eqref{ii19}.
To this end we need also the following lemma,
known as {\em Dirichlet's theorem on diophantine
approximation}; its proof is elementary and can be found 
in \cite[Lemma 2.1]{Vaug97}.
\begin{lemma}[Dirichlet]\label{lii1}
  Let $\alpha$ and $Q$ be real and $Q\ge 1$. There
  exist integers $a$ and $q$ such that
  \[
    1 \le q \le Q, \quad (a, q) = 1, \quad
    |q\alpha - a| < Q^{-1}.
  \]
\end{lemma}

Let $\alpha \in \mathfrak m$. By \eqref{ii4} and Lemma \ref{lii1}
with $Q =  nP^{-1} $, there are integers $a$ and $q$
such that
\[
  P < q \le nP^{-1}, \quad (a, q) = 1, \quad
  |q\alpha - a| < Pn^{-1} \le q^{-1}.
\]
Hence, an appeal to \eqref{ddd1} and Lemma \ref{lii2} gives 
\begin{equation}\label{ii19a}
  f(\alpha) \ll (\log n)^3 \left( nP^{-1/2} + n^{4/5} \right)
    \ll n (\log n)^{ 3 - B/2}.
\end{equation}
and \eqref{ii19} follows on choosing $B\ge 12$. This
completes the proof of Theorem~\ref{th1}.

\medskip

The above proof of Vinogradov's theorem employs the Siegel--Walfisz
theorem and, therefore, is ineffective (recall the
remark following the statement of Theorem \ref{SW}). The interested
reader can find an effective proof (with a slightly weaker error
term) in \cite[Chapter~10]{Kara93}.

\subsubsection{The estimation of $f(\alpha)$}
\label{ssii1.4}

The main tool in the proof of Lemma \ref{lii2} are estimates for
bilinear sums of the form
\begin{equation}\label{ii20}
  S = \mathop{\sum_{X < x \le 2X} \sum_{Y < y \le 2Y}}_{xy \le n}
  \xi_x \eta_y e(\alpha xy).
\end{equation}
We need to control two kinds of such sums, known as \emph{type I
sums} and \emph{type II sums}. For simplicity, we describe these
two types of sums in the simplest cases, noting that the more
general sums arising in the actual proof of Lemma \ref{lii2} can
be reduced to these special cases using standard trickery:
\begin{itemize}
  \item {type I sums:} sums \eqref{ii20} with
       $|\xi_x| \le 1$, $\eta_y = 1$ for all $y$, and
       $X$ is ``not too large'';
  \item {type II sums:} sums \eqref{ii20} with
       $|\xi_x| \le 1$, $|\eta_y| \le 1$, and
       $X, Y$ are ``neither large, nor small''.
\end{itemize}

Vinogradov reduced the estimation of $f(\alpha)$ to the estimation
of type I and type II sums by means of an intricate combinatorial
argument. Nowadays we can achieve the same result almost 
instantaneously by referring to the combinatorial identities of 
Vaughan~\cite{Vaug75, Vaug77} or Heath-Brown~\cite{H-Br82}.
Let $\Lambda(k)$ denote von Mangoldt's function, whose value is 
$\log p$ or $0$ according as $k$ is a power of a prime $p$ or not.
Vaughan's identity states that if $U$ and $V$ are real parameters
exceeding $1$, then
\begin{equation}\label{ii20a}
  \Lambda(k) = \sum_{ \substack{ dm = k\\ 1 \le d \le V}} 
    \mu(d) \log m 
  - \sum_{ \substack{ dlm = k\\ 1 \le d \le V \\ 1 \le m \le U}} 
    \mu(d) \Lambda(m) 
  - \sum_{ \substack{ dlm = k\\ 1 \le d \le V \\ m > U, dl > V}} 
    \mu(d) \Lambda(m).          
\end{equation}
Heath-Brown's identity states that if $k \le x$ and $J$ is a 
positive integer, then
\[
  \Lambda(k) = \sum_{j = 1}^J \binom Jj (-1)^{j - 1}
  \sum_{\substack{ m_1 \cdots m_{2j} = k \\ m_1, \dots, m_j \le x^{1/J} }}
  \mu(m_1) \cdots \mu(m_j) \log m_{2j},
\]
where $\mu(m)$ is the M\"obius function.

Both identities can be used to reduce $f(\alpha)$ to type I
and type II sums with equal success. Here, we apply Vaughan's
identity with $U = V = n^{2/5}$. We obtain
\begin{equation}\label{ii20b}
  \sum_{k \le n} \Lambda(k) e (\alpha k) = W_1 - W_2 - W_3,
\end{equation}
with
\[
  W_j = \sum_{k \le n} a_j(k) e(\alpha k) 
  \quad (1 \le j \le 3)
\]
where $a_j(k)$ denotes the $j$th sum on the right side of
\eqref{ii20a}. The estimation of the sum on the left side of
\eqref{ii20b} is essentially equivalent to that of $f(\alpha)$.
The sums $W_1$ and $W_2$ on the right side of \eqref{ii20b} can
be reduced to type I sums with $X \ll n^{4/5}$; $W_3$ 
can be reduced to type II sums with $n^{2/5} \ll X, Y \ll n^{3/5}$. 
The reader can find all the details in the proof of
\cite[Theorem 3.1]{Vaug97}. Here we will be content with a
brief description of the estimation of the type I and type II
sums.

Consider a type I sum $S_1$. We have
\begin{equation}\label{ii21}
  |S_1| \le \sum_{X < x \le 2X}
  \Big| \sum_{Y < y \le Y'} e (\alpha xy) \Big|,
\end{equation} 
where $Y' = \min (2Y, n/x)$. We can estimate the inner sum
in \eqref{ii21} by means of the elementary bound
\begin{equation}\label{ii22}
  \Big| \sum_{a < y \le b} e (\alpha  y) \Big|
  \le \min \left( b - a + 1, \| \alpha \|^{-1} \right),
\end{equation}
where $\| \alpha \|$ denotes the distance from $\alpha$ to the
nearest integer.  This inequality
follows on noting that the sum on 
the left is the sum of a geometric progression. We obtain 
\begin{equation}\label{ii23}
  |S_1| \le \sum_{x \le 2X} 
  \min \left( Y, \| \alpha x \|^{-1} \right) = T(\alpha),
  \qquad \text{say}.
\end{equation}
Obviously, the trivial estimate for $T(\alpha)$ is
\[
   T(\alpha) \ll XY .
\]
However, under the hypotheses of Lemma~\ref{lii2},
one can establish by elementary methods that (see \cite[Lemma 2.2]{Vaug97})
\begin{equation}\label{ii23a}
  T(\alpha) \ll XY \left( \frac 1q + \frac 1Y + \frac {q}{XY} \right) \log (2XYq).
\end{equation}
Inserting this bound into the right side of \eqref{ii23}, we obtain a 
satisfactory bound for $S_1$.

To estimate a type II sum $S_2$, we first apply Cauchy's 
inequality and get
\[
  |S_2|^2 \ll Y \sum_{Y < y \le 2Y} 
  \Big| \sum_{X < x \le X'} \xi_x e (\alpha xy) \Big|^2,
\]
where $X' = \min (2X, n/y)$. Squaring out and interchanging
the order of summation, we deduce 
\begin{align*}
  |S_2|^2 &\ll Y \sum_{Y < y \le 2Y} \sum_{X < x_1, x_2 \le X'}
    \xi_{x_1} \overline \xi_{x_2} e( \alpha (x_1 - x_2)y ) \\
  &\ll Y \sum_{X < x_1, x_2 \le 2X} \Big| \sum_{Y < y \le Y'} 
    e( \alpha (x_1 - x_2)y ) \Big| \\
  &\ll Y \sum_{X < x \le 2X} \sum_{|h| < X} 
    \Big| \sum_{Y < y \le Y'} e( \alpha hy ) \Big|, 
\end{align*}
where $Y < Y' \le 2Y$. We remark that the innermost sum is now free 
of ``unknown'' weights and can be estimated by means of \eqref{ii22}.
We get
\begin{equation}\label{ii24}
  |S_2|^2 \ll XY^2 + XY T(\alpha),
\end{equation}
and \eqref{ii23a} again leads to a satisfactory bound for $S_2$.

\subsection{The exceptional set in Goldbach's problem}
\label{ssii5}

We now sketch the proof of \eqref{1.5}. We will not discuss the proof
of the more sophisticated results of Montgomery and Vaughan \cite{MoVa75}
and Pintz \cite{Pint04}, since they require knowledge of the properties
of Dirichlet $L$-functions far beyond the scope of this survey. The 
reader can find excellent expositions of the Montgomery--Vaughan result
in their original paper and also in the monograph \cite{PaPa89}.

For an even integer $n$, let $r(n)$ denote the number of
representations of $n$ as the sum of two primes, let $\mathcal 
Z(N)$ denote the set of even integers $n \in (N, 2N]$ with 
$r(n) = 0$, and write $Z(N) = |\mathcal Z(N)|$. Since
\[
  E(x) = \sum_{j = 1}^{\infty} Z(x2^{-j}),
\]
it suffices to bound $Z(N)$ for large $N$.

Define $f(\alpha)$, $\mathfrak M$, and $\mathfrak m$ as before,
with $N$ in place of $n$. When $n$ is an even integer in $(N, 
2N]$, a variant of the method in \S \ref{ssii1.2} gives 
\[
  \int_{\mathfrak M} f(\alpha)^2 e(-\alpha n) \, d\alpha
  = \mathfrak S_{2}(n) \frac n{(\log n)^2} + 
  O\left( \frac N{(\log N)^3} \right),
\]
where 
\[
  \mathfrak S_2(n) = \prod_{p \nmid n} \left( 1 - \frac 1{(p - 1)^2} \right)
  \prod_{ p \mid n } \left( \frac p{p - 1} \right)
\]
is the singular series. In particular, we have $\mathfrak S_{2}(n) \ge 1$ for
even $n$. Thus, for $n \in \mathcal Z(N)$, we have
\begin{equation}\label{ii24a}
  \bigg| \int_{\mathfrak m} f(\alpha)^2 e(-\alpha n) \, d\alpha 
  \bigg| = \bigg| -\int_{\mathfrak M} f(\alpha)^2 e(-\alpha n) 
  \, d\alpha \bigg| \gg N(\log N)^{-2},
\end{equation}
whence
\begin{equation}\label{ii42}
  Z(N) \ll N^{-2}(\log N)^4
  \sum_{n \in \mathcal Z(N)} \bigg| 
  \int_{\mathfrak m} f(\alpha)^2 e(-\alpha n) \, d\alpha \bigg|^2.
\end{equation}
On the other hand, by Bessel's inequality, 
\begin{equation}\label{ii43}
  \sum_{n \in \mathcal Z(N)} \bigg| 
  \int_{\mathfrak m} f(\alpha)^2 e(-\alpha n) \, d\alpha \bigg|^2
  \le \int_{\mathfrak m} |f(\alpha)|^4 \, d\alpha,
\end{equation}
and \eqref{ii18} and \eqref{ii19a} yield
\begin{equation}\label{ii43a}
  \int_{\mathfrak m} |f(\alpha)|^4 \, d\alpha
  \le \Big( \sup_{ \alpha \in \mathfrak m } |f(\alpha)| \Big)^2 
  \int_0^1 |f(\alpha)|^2 \, d\alpha \ll N^3P^{-1}(\log N)^5.
\end{equation}    
Combining \eqref{ii42}--\eqref{ii43a}, we conclude that
\[
  Z(N) \ll NP^{-1}(\log N)^9 \ll N(\log N)^{-A},
\]
on choosing, say, $P = (\log N)^{A + 9}$. This completes the proof
of \eqref{1.5}.

\subsection{The circle method in the Waring--Goldbach problem}
\label{ssii2}

We now turn our attention to Theorems \ref{th3} and \ref{th4}. Much of
the discussion in \S \ref{ssii1} can be generalized to $k$th powers
($k \ge 2$). Using \eqref{ii1}, we can write $R_{k, s}^*(n)$ as
\[
  R_{k, s}^*(n) =
  \int_0^1 f(\alpha)^s e(-\alpha n) \, d\alpha,
\]
where now
\[
  f(\alpha) = \sum_{p \le N} e \left( \alpha p^k \right),
  \quad N = n^{1/k}.
\]
Define the sets of major and minor arcs as before (that is, by 
\eqref{ii3} and \eqref{ii4}, with $P = (\log N)^B$ and $B = B(k, s)$ 
to be chosen later). The machinery in \S \ref{ssii1.2} generalizes 
to $k$th powers with little extra effort. The argument leading to 
\eqref{ii9} gives
\begin{equation}\label{ii26}
  f \left( \frac aq + \beta \right) =
  \phi(q)^{-1}S(q, a) \, v(\beta) \, + \text{ error term},
\end{equation}
where $S(q, a)$ is defined by \eqref{1.13} and
\[
  v(\beta) = \int_2^N \frac {e \left( \beta u^k \right)}
  {\log u} \, du.
\]
We now raise \eqref{ii26} to the $s$th power and integrate the
resulting approximation for $f(\alpha)^s$ over $\mathfrak M$.
Using known estimates for $v(\beta)$ and $S(q, a)$, we find 
that when $s \ge k + 1$,
\begin{equation}\label{ii29}
  \int_{\mathfrak M} f(\alpha)^s e(-\alpha n) \, d\alpha
  = \mathfrak S_{k, s}^*(n) J_{k, s}^*(n)
  + O \left( N^{s - k} P^{-1/k + \eps} \right),
\end{equation}
where $\mathfrak S_{k, s}^*(n)$ is defined by \eqref{1.13} and
$J_{k, s}^*(n)$ is the singular integral
\begin{align*}
    J_{k, s}^*(n) 
  &= \int_{-\infty}^{\infty}  v(\beta)^s e(-\beta n) \, d\beta \\
  &= \frac {\Gamma^s \left( 1 + {\textstyle \frac 1k} \right)}{\Gamma\left(
\textstyle \frac sk \right)}
    \frac {n^{s/k - 1}}{(\log n)^s} +
    O\left( n^{s/k - 1}(\log n)^{-s - 1} \right).
\end{align*}
This reduces the proof of Theorem \ref{th3} to the estimate
\begin{equation}\label{ii31}
  \int_{\mathfrak m} f(\alpha)^s e(-\alpha n) \, d\alpha
  \ll N^{s - k} (\log N)^{-s - 1}.
\end{equation}

Notice that when $k = 1$ and $s = 3$, \eqref{ii31} turns into
\eqref{ii16}. Thus, it is natural to try to obtain variants of 
\eqref{ii18} and \eqref{ii19} for $f(\alpha)$ when $k \ge 2$. 
To estimate the maximum of $f(\alpha)$ on the minor arcs, we use 
the same tools as in \S \ref{ssii1.4}, that is:
\begin{itemize}
  \item Heath-Brown's or Vaughan's identity to reduce the estimation
       of $f(\alpha)$ to the estimation of bilinear sums
       \[
         \mathop{\sum_{X < x \le 2X} \sum_{Y < y \le 2Y}}_{xy \le N}
         \xi_x \eta_y e \left( \alpha (xy)^k \right);
       \]
  \item Cauchy's inequality to bound those bilinear sums  
       in terms of the quantity $T(\alpha)$ appearing in \eqref{ii23}.
\end{itemize} 
The following result due to Harman~\cite{Harm81} is the analogue of
Lemma \ref{lii2} for $k \ge 2$.

\begin{lemma}\label{lii3}
  Let $k \ge 2$, let $\alpha \in \mathbb R$, and suppose that $a$
  and $q$ are integers satisfying
  \[
    1 \le q \le N^k, \quad (a, q) = 1, \quad
    |q\alpha - a| < q^{-1}.
  \]
  There is a constant $c = c(k) > 0$ such that
  \[
    f(\alpha) \ll N(\log N)^c
    \left( q^{-1} + N^{-1/2} + qN^{-k} \right)^{4^{1 - k}}.
  \]
\end{lemma}

On choosing the constant $B$ (in the definition of $\mathfrak m$)
sufficiently large, one can use Lemmas \ref{lii1} and \ref{lii3} 
to show that, for any fixed $A > 0$,
\[
  \sup_{\alpha \in \mathfrak m} |f(\alpha)|
  \ll N (\log N)^{-A}.
\]
Hence, if $s = 2r + 1$, one has
\[
  \int_{\mathfrak m} |f(\alpha)|^s \, d\alpha
  \le \sup_{\alpha \in \mathfrak m} |f(\alpha)| \int_0^1 |f(\alpha)|^{2r} \, d\alpha
  \ll N(\log N)^{-A} \int_0^1 |f(\alpha)|^{2r} \, d\alpha,
\]
and it suffices to establish the estimate
\begin{equation}\label{ii33}
  I_r(N) := \int_0^1 |f(\alpha)|^{2r} d\alpha
  \ll N^{2r - k} (\log N)^c,
\end{equation}
with $c = c(k, r)$.

\subsection{Mean-value estimates for exponential sums}
\label{s3.4}

We now turn to the proof of \eqref{ii33}. By \eqref{ii1}, 
$I_r(N)$ represents the number of solutions of the diophantine
equation
\begin{equation}\label{ii35}
  \begin{cases}
    x_1^k + \cdots + x_r^k = x_{r + 1}^k + \cdots + x_{2r}^k, \\
    1 \le x_1, \dots, x_{2r} \le N
  \end{cases}
\end{equation}
in {\em primes} $x_1, \dots, x_{2r}$, and therefore, $I_r(N)$
does not exceed the number of solutions of \eqref{ii35} in
{\em integers} $x_1, \dots, x_{2r}$. Using \eqref{ii1} to write
the latter quantity as a Fourier integral, we conclude that
\begin{equation}\label{ii36}
  I_r(N) \le \int_0^1 |g(\alpha)|^{2r} d\alpha, \qquad
  g(\alpha) = \sum_{x \le N} e \left( \alpha x^k \right).
\end{equation}
This reduces the estimation of the even moments of $f(\alpha)$
to the estimation of the respective moments of the exponential
sum $g(\alpha)$, whose analysis is much easier. In particular,
we have the following two results.

\begin{lemma}[Hua's lemma]\label{lii4}
  Suppose that $k \ge 1$, and let $g(\alpha)$ be defined by
  \eqref{ii36}. There exists a constant $c = c(k) \ge 0$
  such that
  \begin{equation}\label{3.39}
    \int_0^1 |g(\alpha)|^{2^k} d\alpha \ll
    N^{2^k - k}(\log N)^c.
  \end{equation}
\end{lemma}

\begin{lemma}\label{lii5}
  Suppose that $k \ge 11$ and $g(\alpha)$ is defined by \eqref{ii36}. 
  There exists a constant $c = c(k) > 0$ such that for $r > \frac 12
  k^2( \log k + \log\log k + c)$, 
  \begin{equation}\label{3.40}
    \int_0^1 |g(\alpha)|^{2r} d\alpha \ll N^{2r - k}.
  \end{equation}
\end{lemma}

These lemmas are, in fact, rather deep and important results in
the theory of Waring's problem. Unfortunately, their proofs are
too complicated to include in this survey in any meaningful way.
The reader will find a proof of a somewhat weaker version of Hua's
lemma (with a factor of $N^{\eps}$ in place of $(\log N)^c$) in
\cite[Lemma 2.5]{Vaug97} and a complete proof in \cite[Theorem 4]
{Hua65}. Results somewhat weaker than Lemma \ref{lii5} are 
classical and go back to Vinogradov's work on Waring's problem
(see \cite[Lemma 7.13]{Hua65} or \cite[Theorem 7.4]{Vaug97}). 
Lemma \ref{lii5} itself follows from the results in Ford 
\cite{Ford95} (in particular, see \cite[(5.4)]{Ford95}).

Combining \eqref{ii36} and Lemmas \ref{lii4} and \ref{lii5}, we get 
\eqref{ii33} with 
\[
  r = \begin{cases}
    2^{k - 1}, & \text{ if } \;\; k \le 10, \\
    \left[ \frac 12k^2 (\log k + \log \log k + c) \right] + 1,
               & \text{ if } \;\; k \ge 11.
  \end{cases}
\]
Clearly, this completes the proof of Theorem \ref{th3}, except 
for the case $6 \le k \le 8$, which we will skip in order to
avoid the discussion of certain technical details.

\subsection{Diminishing ranges}
\label{ssii4}

In this section, we describe the main new idea that leads to
the bounds for $H(k)$ in Theorem \ref{th4}. This idea, known as
the \emph{method of diminishing ranges}, appeared for the first
time in the work of Hardy and Littlewood on Waring's problem
and later was developed into a powerfull technique by Davenport.

The limit of the method employed in \S\ref{ssii2} is set by the
mean-value estimates  in Lemmas \ref{lii4} and \ref{lii5}.
The key observation
in the method of diminishing ranges is that it can be much easier
to count the solutions of the equation in \eqref{ii35} if the
unknowns $x_1, \dots, x_{2r}$ are restricted to proper subsets
of $[1, N]$. For example, the simplest version of the method that
goes back to Hardy and Littlewood uses that when $N_2, \dots,
N_r$ are defined recursively by
\[
  N_j = k^{-1}N_{j - 1}^{1 - 1/k} \qquad (2 \le j \le r),
\]
the equation
\begin{equation}\label{ii35a}
  \begin{cases}
    x_1^k + \cdots + x_r^k = x_{r + 1}^k + \cdots + x_{2r}^k, \\
    N_j < x_j, x_{r + j} \le 2N_j \quad (1 \le j \le r),
  \end{cases}
  \tag{\ref{ii35}*}
\end{equation}
has only ``diagonal" solutions with $x_{r + j} = x_j$, $j = 1,
\dots, r$. Thus, the number of solutions of \eqref{ii35a} is
bounded above by
\[
  N_1 \cdots N_r \ll N_1^{2-\lambda}(N_2 \cdots N_r)^2
\]
where
\[
  \lambda = 1 + {\textstyle \big( 1 - \frac 1k \big) +
  \dots + \big( 1 - \frac 1k \big)^{r - 1}} \ge k - ke^{-r/k}.
\]
That is, we have the bound
\begin{equation}\label{3.41}
  \int_0^1
  \big| g_1(\alpha)g_2(\alpha) \cdots g_r(\alpha) \big|^2
  d\alpha \ll N_1^{2 - \lambda}(N_2 \cdots N_r)^2,
\end{equation}
where
\[
  g_j(\alpha) = \sum_{N_j < x \le 2N_j}
  e \left( \alpha x^k \right) \quad (1 \le j \le r).
\]

We can use \eqref{3.41} as a replacement
for the mean-value estimates in \S \ref{s3.4}. Let $T_{k, s}(n)$
denote the number of solutions of
\[
  p_1^k + p_2^k + \dots + p_s^k = n
\]
in primes $p_1, \dots, p_s$ subject to
\[
  N_j < p_j, p_{r + j} \le 2N_j \quad (1 \le j \le r), \qquad
  N_1 < p_{2r + 1}, \dots, p_s \le 2N_1.
\]
Then
\begin{equation}\label{3.42}
  T_{k, r}(n) = \int_0^1
  f_1(\alpha)^{s - 2r + 2} f_2(\alpha)^2 \cdots f_r(\alpha)^2
  e( -\alpha n ) \, d\alpha,
\end{equation}
where
\[
  f_j(\alpha) = \sum_{N_j < p \le 2N_j} e \left( \alpha p^k \right)
  \quad (1 \le j \le r).
\]
When $r \sim ck\log k$, we can use \eqref{3.41} to derive a bound
of the form
\[
  \int_0^1
  \big| f_1(\alpha)^2 f_2(\alpha) \cdots f_r(\alpha) \big|^2
  d\alpha \ll N_1^{4 - k} (N_2 \cdots N_r)^2.
\]
Furthermore, assuming that $s$ is just slightly larger than $2r$
(it suffices to assume that $s \ge 2r + 3$, for example), we can
then obtain an asymptotic formula for the right side of
\eqref{3.42} by the methods sketched in \S \ref{ssii2}. This is
(essentially) how one proves Theorem \ref{th4} for $k \ge 11$.
The proof for $k \le 10$ follows the same general approach, except
that we use more elaborate choices of the parameters $N_1, \dots,
N_r$ in \eqref{ii35a}.

\subsection{Kloosterman's refinement of the circle method}

Consider again equation \eqref{1.9} with $k=2$. The Hardy--Littlewood 
method in its original form establishes the asymptotic formula
\eqref{1.10} for $s > 4$, but it fails to prove Lagrange's four 
squares theorem. In 1926
Kloosterman~\cite{Kloo26} proposed a variant of the circle
method, known today as {\it Kloosterman's refinement}, which he
used to prove an asymptotic formula for the number of solutions
of the equation
\begin{equation} \label{d10}
  a_1 x_1^2 + \dots + a_4 x_4^2 = n ,
\end{equation}
where $a_i$ are fixed positive integers.

Denote by $I(n)$ the number of solutions
of \eqref{d10} in positive integers $x_i$.
By \eqref{ii1},
\begin{equation} \label{d11}
  I(n) = \int_0^1 H(\alpha)e(-\alpha n) \, d \alpha,
\end{equation}
where
\[
 H(\alpha) = h(a_1 \alpha) \cdots h(a_4 \alpha),
 \qquad
 h(\alpha) = \sum_{ x \le N } e\left(\alpha x^2\right) ,
 \qquad
 N = n^{1/2}.
\]
A ``classical'' Hardy--Littlewood decomposition of the right
side of \eqref{d11} into integrals over major and minor arcs
is of little use here, since we cannot prove that the
contribution from the minor arcs is smaller than the expected
main term. Kloosterman's idea is to eliminate the minor arcs
altogether.

The elimination of the minor arcs requires greater care in
the handling of the major arcs. Let $X$ be the integer with
$X - 1 < N \le X$. It is clear that the integration in 
\eqref{d11} can be taken over the interval $\big( X^{-1}, 
1 + X^{-1} \big]$, which can be represented as a union of
disjoint intervals
\begin{equation} \label{d12.5}
  \left( X^{-1}, 1 + X^{-1} \right] =
  \bigcup_{q \le N}
  \bigcup_{ \substack{ 1 \le a \le q\\ (a, q) = 1}}
  \left( \frac aq - \frac 1{qq_1}, \frac aq + \frac 1{qq_2} \right],
\end{equation}
where for each pair $q, a$ in the union, the positive integers
$q_1 = q_1(q, a)$ and $q_2 = q_2(q, a)$ are uniquely determined 
and satisfy the conditions
\begin{equation} \label{d11.5}
  N < q_1, q_2 \le 2N, \qquad aq_1 \equiv 1 \pmodulo q, \qquad
  aq_2 \equiv -1 \pmodulo q.
\end{equation}
The decomposition \eqref{d12.5} is known as the \emph{Farey 
decomposition} and provides a natural way of partitioning of
the unit interval into non-overlapping major arcs (see Hardy and 
Wright \cite[Section 3.8]{HaWr79}). Let $\mathfrak M(q, a)$ 
denote the interval in the Farey decomposition ``centered'' at
$a/q$. We have
\begin{equation} \label{d12}
 I(n) = \sum_{ q \le N}
 \sum_{ \substack{ 1 \le a \le q \\ (a, q)=1 }}
 e \Big( \frac {-an}q \Big)
 \int_{ \mathfrak B (q, a) }
  H \Big( \frac{a}{q} + \beta \Big) e( -\beta n) \, d\beta ,
\end{equation}
where $\mathfrak B(q, a)$ is defined by
\begin{equation} \label{d13}
  \mathfrak B(q, a) =
  \left\{ \beta \in \mathbb R : 
  a/q + \beta \in \mathfrak M(q, a) \right\}.
\end{equation}

We can find an asymptotic formula for the integrand on the 
right side of \eqref{d12}. The contribution of the main term 
in that asymptotic formula produces the expected main term 
in the asymptotic formula for $I(n)$. However, in order to 
obtain a satisfactory bound for the contribution of the 
error term, we have to take into account the cancellation 
among terms corresponding to different Farey fractions $a/q$
with the same denominator. To this end, we want to interchange
the order of integration and summation over $a$ in \eqref{d12}.
Since the endpoints of $\mathfrak B(q, a)$ depend on $a$, the
total contribution of the error terms can be expressed as
\begin{equation} \label{a15}
  \sum_{q \le N} \int_{-1/(qN)}^{1/(qN)} 
  \bigg\{\; \sideset{}{^{(\beta)}}
  \sum_{ \substack{ 1 \le a \le q\\ (a, q) = 1}}
  E \left( \frac aq + \beta \right) e \left( \frac {-an}q \right) \bigg\}
  e(-\beta n) \, d\beta,
\end{equation}
where $E(a/q + \beta)$ is the error term in the major arc
approximation to $H(a/q + \beta)$ and the superscript in 
$\sum^{(\beta)}$ indicates that the summation is restricted 
to those $a$ for which $\mathfrak B (q, a) \ni \beta$. Using 
\eqref{d11.5} and \eqref{d13}, we can transform the latter 
constraint on $a$ into a condition about the multiplicative 
inverse of $a$ modulo $q$, that is, the unique residue class 
$\bar a$ modulo $q$ with $\bar aa \equiv 1 \pmodulo q$. Thus, 
a special kind of exponential sums enter the scene: the 
\emph{Kloosterman sums}
\[
  K(q; m, n) = \sum_{ \substack{ x=1 \\ (x,q)=1 }}^q
  e \Big( \frac{ m x + n \bar x }{q} \Big) .
\]
There also other (in fact, more substantial) reasons for the 
Kloosterman sums to appear, but those are too technical to
include here.

The success of Kloosterman's method hinges on the existence of 
sufficiently sharp estimates for $K(q; m, n)$. The first such estimate
was found by Kloosterman himself and later his result
has been improved. Today it is known that 
\begin{equation}\label{d16}
  | K(q; m, n) | \le \tau(q) \, q^{1/2} \, (m, n, q)^{1/2} ,
\end{equation}
where $(m, n, q)$ is the greatest common divisor of $m, n, q$
and $\tau(q)$ is the number of positive divisors of $q$.
In 1948 A. Weil \cite{Weil48} proved \eqref{d16} in the most 
important case: when $q$ is a prime. In the general case \eqref{d16}
was established by Estermann~\cite{Este61}. This estimate plays 
an important role not only in the Kloosterman refinement of the 
circle method, but in many other problems in number theory.

Kloosterman's method has been applied to several additive problems, 
and in particular, to problems with primes and almost primes. We 
refer the reader, for example, to Estermann~\cite{Este62},
Hooley~\cite{Hoo86}, 
Heath-Brown~\cite{H-Br83, H-Br96, H-Br98},
Br\"udern and Fouvry~\cite{BrFo94},
Heath-Brown and Tolev~\cite{H-BTo03}.

\section{Sieve methods}
\label{siv}

In this section we describe the so-called \emph{sieve methods},
which are an important tool in analytic number theory and, in 
particular, in the proof of Chen's theorem (Theorem \ref{th2}
in the Introduction). We start with a brief account of the main
idea of the method (\S \ref{siv0} and \S \ref{siv1}). This allows 
us in \S \ref{siv2} to present a proof of a slightly weaker (but 
much simpler) version of Chen's result, in which $P_2$-numbers 
are replaced by $P_4$-numbers. We conclude the section by 
sketching some of the new ideas needed to obtain Chen's theorem 
in its full strength (\S \ref{siv3}) and of some further work on 
sieve methods (\S \ref{siv4}).

\subsection{The sieve of Eratosthenes}
\label{siv0}

Let $\mathcal A$ be a finite integer sequence. We will be
concerned with the existence of elements of $\mathcal A$ that
are primes or, more generally, almost primes $P_r$, with $r$
bounded. In general, sieve methods reduce such a question to
counting the elements $a \in \mathcal A$ not divisible by
small primes $p$ from some suitably chosen set of primes
$\mathfrak P$. To be more explicit, we consider
a set of prime numbers $\mathfrak P$ and a real parameter
$z \ge 2$ and define the \emph{sifting function}
\begin{equation}\label{iv1}
  S(\mathcal A, \mathfrak P, z) =
  \left| \{ a \in \mathcal A : (a, P(z))=1 \} \right|, \quad
  P(z) = \prod_{\substack{ p < z \\ p \in \mathfrak P}} p,
\end{equation}
where $|\mathcal A|$ denotes the number of elements of a
sequence $\mathcal A$ (not the cardinality of the
underlying set). In applications, the set $\mathfrak P$ is
usually taken to be the set of possible prime divisors of
the elements of $\mathcal A$, so the sifting function
\eqref{iv1} counts the elements of $\mathcal A$ free of
prime divisors $p < z$.

In order to bound $S(\mathcal A, \mathfrak P, z)$, we recall
the following fundamental property of the M\"obius function
(see \cite[Theorem 263]{HaWr79}):
\begin{equation}\label{iv2}
  \sum_{ d \mid k} \mu(d) =
  \begin{cases}
         1,       & \text{if } \;\; k = 1,     \\
         0,       & \text{if } \;\; k > 1.
  \end{cases}
\end{equation}
Using this identity, we can express the sifting function in the form
\begin{equation}\label{iv3}
  S(\mathcal A, \mathfrak P, z) =
  \sum_{a \in \mathcal A} \sum_{d \mid (a, P(z))} \mu(d).
\end{equation}
We can now interchange the order of summation to get
\begin{equation}\label{iv4}
  S(\mathcal A, \mathfrak P, z) =
  \sum_{d \mid P(z)} \mu(d) |\mathcal A_d|,
\end{equation}
where
\[
  \mathcal A_d = \{ a \in \mathcal A : a \equiv 0 \pmodulo d \}.
\]
To this end, we suppose that there exist a (large) parameter $X$
and a multiplicative function $\omega(d)$ such that
$|\mathcal A_d|$ can be approximated by
$ X \omega(d)/ d$.
We write $r(X, d)$ for the error term in
this approximation, that is,
\begin{equation}\label{iv5}
  |\mathcal A_d| = X \frac {\omega(d)}{d} + r(X, d).
\end{equation}
We expect $r(X, d)$ to be `small', at least in some
average sense over $d$. Substituting \eqref{iv5} into the right
side of \eqref{iv4}, we find that
\begin{equation}\label{iv6}
  S(\mathcal A, \mathfrak P, z) = XV(z) + R(X, z),
\end{equation}
where
\begin{equation}\label{iv7}
  V(z) = \sum_{d \mid P(z)} \mu(d) \frac {\omega(d)}{d},
  \quad  R(X, z) = \sum_{d \mid P(z)} \mu(d) r(X, d).
\end{equation}
We would like to believe that, under `ideal circumstances',
\eqref{iv6} is an asymptotic formula for the sifting function
$S(\mathcal A, \mathfrak P, z)$, $XV(z)$ being the main
term and $R(X, z)$ the error term. However, such expectations
turn out to be unrealistic, as we are about to demonstrate.

Let us try to apply \eqref{iv6} to bound above
the number of primes $\le x$. We choose
\begin{equation}\label{iv8}
  \mathcal A = \{ n \in \mathbb N : n \le x \}, \quad
  \mathfrak P = \{p: p \text{ is a prime} \}.
\end{equation}
Then
\[
  |\mathcal A_d| = \left[ \frac xd \right] =
  \frac xd + r(x, d), \quad |r(x, d)| \le 1,
\]
that is, $X = x$ and $\omega(d) = 1$ for all $d$.
Using an elementary property of multiplicative functions
(see \cite[Theorem 286]{HaWr79}), we can write $V(z)$ as
\begin{equation}\label{iv9}
  V(z) = \prod_{p < z} \left( 1 - \frac{\omega(p)}{p} \right).
\end{equation}
When $\omega(p) = 1$, this identity and an asymptotic formula
due to Mertens (see \cite[Theorem 429]{HaWr79})
reveal that the main term in \eqref{iv6} is
\begin{equation}\label{iv10}
  XV(z) = X\prod_{p < z} \left( 1 - \frac 1p \right)
  \sim X\frac {e^{-\gamma}}{\log z} \qquad \text{as } \quad z \to \infty;
\end{equation}
here $\gamma = 0.5772\dots$ is Euler's constant. Thus,
if $\mathcal A$ and $\mathfrak P$ are as in \eqref{iv8} and
$z = x^{1/2}$, the projected `main term' in \eqref{iv6} is
$\sim 2e^{-\gamma}x(\log x)^{-1}$ as $x \to \infty$, whereas
the true size of the sifting function on the left side is
\[
  S(\mathcal A, \mathfrak P, \sqrt x) = \pi(x) - \pi(\sqrt x) + 1
  \sim \frac x{\log x}  \qquad \text{as } \quad x \to \infty,
\]
by the Prime Number Theorem. Since $2e^{-\gamma} = 1.122\dots$,
we conclude that the `error term' $R(x, \sqrt x)$ is in this case
of the same order of magnitude as the `main term'.

Identity \eqref{iv6} is known as the \emph{sieve of
Eratosthenes--Legendre}. The basic idea goes back to the ancient Greeks
(usually attributed to Eratosthenes), while the formal exposition above
is essentially due to Legendre, who used the above argument to show that
\[
  \pi(x) \ll \frac x{\log\log x} .
\]
The sieve of Eratosthenes--Legendre can be extremely powerful in certain
situations\footnote{For example, I.~M.~Vinogradov's combinatorial
argument for converting sums over primes into linear combinations
of type~I and type II sums is based on a variant of \eqref{iv6}.
See Harman~\cite{Harm97} for other applications and further
discussion.}, but in most cases the sum $R(X, z)$ contains `too
many' terms for \eqref{iv6} to be of any practical use (e.g., in
the above example, $R(X, z)$ contains $2^{\pi(z)}$ terms). Modern
sieve methods use various clever approximations to the left side
of \eqref{iv2} to overcome this problem. In the following sections,
we describe one of the variants of one the existing approaches.
The reader can find other constructions, comparisons of the various
approaches, and proofs in the monographs on sieve methods
\cite{Grea01, HaRi74, Moto83} or in \cite{H-Br03}
(see also the remarks in \S\ref{siv4} for other references).

\subsection{The linear sieve}
\label{siv1}

Let $y > 0$ be a parameter to be chosen later in terms of $X$ and
suppose that $\lambda^+(d)$ and $\lambda^-(d)$ are real-valued 
functions supported on the squarefree integers $d$ (i.e.,
$\lambda^{\pm}(d) = 0$ if $d$ is divisible by the square of a prime).
Furthermore, suppose that
\begin{equation} \label{iv11}
  |\lambda^{\pm}(d)|\le 1  \quad \text{and} \quad
  \lambda^{\pm}(d) = 0 \quad \text{for } \; d \ge y,
\end{equation} 
and that
\begin{equation} \label{iv12}
  \sum_{ d \mid n } \lambda^-(d) \le 
  \sum_{ d \mid n } \mu(d)  \le 
  \sum_{ d \mid n } \lambda^+(d) \quad \text{for all } \; n = 1, 2, \dots.
\end{equation}
Using \eqref{iv3} and the left inequality in \eqref{iv12}, we obtain
\[
  S(\mathcal A, \mathfrak P, z)  \ge  
  \sum_{ a \in \mathcal A } \sum_{ d \mid (a, P(z))} \lambda^-(d).
\]
We can interchange the order of summation in the right side of this
inequality and apply \eqref{iv5} and \eqref{iv11} to get the bound
\begin{align*}
  S(\mathcal A, \mathfrak P, z) &\ge
     \sum_{ d \mid P(z) } \lambda^-(d) |\mathcal A_d| = 
     \sum_{ d \mid P(z) } \lambda^-(d)  
     \left( X \frac{\omega(d)}{d} + r(X, d)  \right) \\
  &= X \sum_{ d \mid P(z) } \lambda^-(d) \frac {\omega(d)}d  
     + \sum_{ d \mid P(z) } \lambda^-(d) r(X, d) \ge 
     X \mathcal M^- - \mathcal R,
\end{align*}
where
\begin{equation} \label{iv13}
  \mathcal M^{\pm} = \sum_{ d \mid P(z) } \lambda^{\pm}(d) 
  \frac {\omega (d)}d, \quad
  \mathcal R = \sum_{ \substack{ d \mid P(z)\\ d < y}} |r(X, d)|.
\end{equation}
In a similar fashion, we can use the right inequality in
\eqref{iv12} to estimate the sifting function from above.
That is, we have
\begin{equation} \label{iv14}
  X\mathcal M^- - \mathcal R \le 
  S(\mathcal A, \mathfrak P, z) \le
  X \mathcal M^+ + \mathcal R.
\end{equation}

We are now in a position to overcome the difficulty caused
by the ``error term'' in the Eratosthenes--Legendre sieve.
The sum $\mathcal R$ is similar to the error term $R(X, z)$ 
defined in \eqref{iv7}, but unlike $R(X, z)$ we can use
the parameter $y$ to control the number of terms in 
$\mathcal R$. Thus, our general strategy will be to 
construct functions $\lambda^{\pm}(d)$ which satisfy
\eqref{iv11} and \eqref{iv12} and for which the sums 
$\mathcal M^{\pm}$ are of the same order as the sum 
$V(z)$ defined in \eqref{iv7}. There are various 
constructions of such functions $\lambda^{\pm}(d)$. However,
since it is not our goal to give a detailed treatment of 
sieve theory here, we will simply state one of the modern
sieves in a form suitable for an application to the binary
Goldbach problem.

The sieve method we will use is known as the \emph{Rosser--Iwaniec
sieve}. Its idea appeared for the first time in an unpublished
manuscript by Rosser. The full-fledged version of this sieve was
developed independently by Iwaniec~\cite{Iwan80a, Iwan80b}. 
Suppose that the multiplicative function $\omega$ in 
\eqref{iv5} satisfies the condition
\begin{equation}\label{iv15}
  \prod_{ w_1 \le p < w_2 } \left( 1 - 
  \frac{\omega(p)}{p} \right)^{-1} \le 
  \left( \frac { \log w_2 }{ \log w_1 } \right)^{\kappa}
  \left( 1 + \frac{K}{\log w_1} \right) \quad (2 \le w_1 < w_2),
\end{equation}
where $\kappa > 0$ is an absolute constant known as the \emph{sieve 
dimension} and $K > 0$ is independent of $w_1$ and $w_2$. This 
inequality is usually interpreted as an average bound for the
values taken by $\omega(p)$ when $p$ is prime, since it is consistent
with the inequality $\omega(p) \le \kappa$. In our application of
the sieve to Goldbach's problem, we will have to deal with a 
sequence $\mathcal A$ (given by \eqref{1.7}) for which \eqref{iv15}
holds with $\kappa = 1$, so we will state the Rosser--Iwaniec sieve 
in this special case, in which it is known as the \emph{linear sieve}.

Suppose that $\omega(p)$ satisfies \eqref{iv15} with $\kappa = 1$ and
that
\begin{equation} \label{iv16}
  0 < \omega(p) < p \quad \text{when } p \in \mathfrak P
  \qquad \text{and} \qquad
  \omega(p) = 0 \quad \text{when } p \not\in \mathfrak P.
\end{equation}
We put $\lambda^{\pm}(1) = 1$ and $\lambda^{\pm}(d) = 0$
if $d$ is not squarefree. If $d > 1$ is squarefree and has
prime decomposition $d = p_1 \cdots p_r$, $p_1 > p_2 > \cdots > p_r$,
we define
\begin{align} 
  \lambda^+(d) &= 
  \begin{cases}
    (-1)^r  & \text{if } p_1 \cdots p_{2l} p_{2l + 1}^3 < y 
              \text{ whenever }  0 \le l \le (r - 1)/2,     \\
    0       & \text{otherwise},
  \end{cases}   \label{iv17}  \\
  \lambda^-(d) & = 
  \begin{cases} 
    (-1)^r  & \text{if } p_1 \cdots p_{2l -1 } p_{2l}^3 < y 
              \text{ whenever }  1 \le l \le r/2, \\
    0       & \text{otherwise}. 
  \end{cases}    \label{iv18} 
\end{align} 
It can be shown (see \cite{Grea01, Iwan80a}) that these two 
functions satisfy conditions \eqref{iv11} and \eqref{iv12}. Furthermore,
if the quantities $\mathcal M^{\pm}$ are defined by \eqref{iv13} with
$\lambda^{\pm}(d)$ given by \eqref{iv17} and \eqref{iv18}, we have
\begin{align} 
   V(z) & \le \mathcal M^+ \le V(z) 
   \left( F(s) + O \big( e^{-s} (\log y)^{-1/3} \big) \right) \;\;
      & \text{for }  \;\; s \ge 1,    \label{iv19}        \\ 
   V(z) & \ge \mathcal M^- \ge V(z) 
   \left( f(s) + O \big( e^{-s} (\log y)^{-1/3} \big) \right) \;\;  
      & \text{for }  \;\; s \ge 2,    \label{iv20}
\end{align} 
where $s = \log y / \log z$ and the functions $f(s)$ and $F(s)$ are 
the continuous solutions of 
a system of differential delay equations
(see \cite{Grea01, Iwan80a}).
The analysis of that system reveals that the function $F(s)$ 
is strictly decreasing for $s > 0$, that the function $f(s)$ is 
strictly increasing for $s > 2$, and that
\begin{equation}\label{iv24}
  0 < f(s) < 1 < F(s) \qquad \text{for } \; s > 2.
\end{equation}
Furthermore, both functions are very close to $1$ for large $s$.
More precisely, they satisfy
\begin{equation} \label{iv25}
  F(s), \; f(s) = 1 + O(s^{-s}) \qquad \text{as } \;\;
  s \to \infty.
\end{equation}
Substituting \eqref{iv19}
and \eqref{iv20} into \eqref{iv14}, we obtain
\begin{align} 
  S(\mathcal A, \mathfrak P, z) 
       &\le
  X V(z) \left( F(s) + O \big( ( \log y )^{-1/3} \big) \right) 
  + \mathcal R \;\;
       & \text{for } \;\;  s \ge 1, \label{iv22} \\
  S(\mathcal A, \mathfrak P, z) 
       & \ge
  X V(z) \left( f(s) + O \big( ( \log y )^{-1/3} \big) \right) 
  - \mathcal R \;\; 
       & \text{for } \;\; s \ge 2, \label{iv23}
\end{align} 
where $\mathcal R$ is defined by \eqref{iv13}.

We now return to our initial goal---namely, to prove that the
sequence $\mathcal A$ contains almost primes. We want to use
\eqref{iv23} to show that
\begin{equation} \label{iv25a}
  S(\mathcal A, \mathfrak P, X^{\alpha}) > 0
\end{equation}
for some fixed $\alpha > 0$. This will imply the existence of 
an $a \in \mathcal A$ all of whose prime divisors exceed $X^{\alpha}$.
If $|a| \ll X^g$ for all $a \in \mathcal A$, it will then follow that 
$\mathcal A$ contains a $P_r$-number, where $r \le g/\alpha$. Clearly, 
since we want to minimize $r$, we would like to take $\alpha$ as large 
as possible. On the one hand, in order to derive \eqref{iv25a} from 
\eqref{iv23}, we need to ensure that the main term in \eqref{iv23} is 
positive and that the error term $\mathcal R$ is of a smaller order of 
magnitude than the main term. It is the balancing of these two 
requirements that determines the optimal choice for $z$ and, ultimately,
the quality of our result. In view of \eqref{iv24}, the positivity of
the main term in \eqref{iv23} requires choosing $y$ slightly larger 
than $z^2$. On the other hand, while in some applications the estimation
of $\mathcal R$ is easier than in others, it is always the case that it
imposes a restriction on how large we can choose $y$, and hence, how 
large we can choose $z$. In the next section, we demonstrate how this
general approach works when applied to the binary Goldbach problem.

\subsection{The linear sieve in the binary Goldbach problem}
\label{siv2}

In this section, we apply the linear Rosser--Iwaniec sieve to
the sequence $\mathcal A$ in \eqref{1.7} and the set $\mathfrak
P$ of odd primes that do not divide $n$, that is,
\[
  \mathcal A = \mathcal A(n) 
  = \left\{ n - p : 2 < p < n \right\} \qquad \text{and} \qquad  
  \mathfrak P = \{ p : p > 2, \; p \nmid n \}.
\]
It is clear that all elements of $\mathcal A$ are odd numbers 
and that at most $\log n$ of them may have a common prime factor 
with $n$ (for $(n, n - p) > 1$ implies $p \mid n$, and $n$ has at 
most $\log n$ odd prime factors). Thus, $\mathfrak P$ is the set
of ``typical'' prime divisors of elements of $\mathcal A$.

Next, we proceed to define the quantity $X$ and the multiplicative 
function $\omega(d)$ in \eqref{iv5}. We have
\begin{equation} \label{iv26}
  |\mathcal A_d| = 
  \sum_{ \substack{ 2 < p < n \\ p \equiv n \pmodulo d}} 1
  = \pi(n; d, n) - 1,
\end{equation}
so the prime number theorem for arithmetic progressions suggests
the choice
\begin{equation} \label{iv27}
  X = \li n \quad \text{and} \quad
  \omega(d) = \begin{cases}
    d / \phi(d)  & \text{if } \; (d, n) = 1,\\
    0            & \text{otherwise}.
  \end{cases}
\end{equation}
With this choice, the error terms $r(X, d)$ defined by \eqref{iv5}
satisfy the inequality
\[
  |r(X, d)| \le \begin{cases}
     \displaystyle 1 + 
     \left| \pi(n; d, n) - \frac {\li n}{\phi(d)} \right| 
     & \text{if } \; (d, n) = 1, \\
     1 & \text{otherwise}.
  \end{cases}
\]
It then follows from the Bombieri--Vinogradov theorem (Theorem
\ref{B-V}) that
\begin{equation}\label{iv28}
  \mathcal R \le y + \sum_{d \le y} \max_{(a, d) = 1} 
  \left| \pi(n; d, a) - \frac {\li n}{\phi(d)} \right|
  \ll n(\log n)^{-3},
\end{equation}
whenever $y \le n^{1/2}(\log n)^{-6}$. Furthermore, we have
\begin{equation}\label{iv29}
  V(z) = \prod_{ \substack{ p < z \\ p \nmid n }}
  \left( 1 - \frac 1{p - 1} \right) 
  \ge \prod_{ p < z } \left( 1 - \frac{1}{p-1} \right)
  \gg (\log z)^{-1}.
\end{equation}
On choosing $y = n^{1/2}(\log n)^{-6}$ and $z = n^{2/9}$, we have
\[
  \frac{ \log y }{\log z } = \frac 94 
  + O \left( \frac{ \log \log n }{\log n} \right) > 2.2,
\]
provided that $n$ is sufficiently large. Hence, we deduce from 
\eqref{iv24}, \eqref{iv23} and \eqref{iv27}--\eqref{iv29} that
\begin{equation} \label{sieve}
  S(\mathcal A, \mathfrak P, z) \gg n( \log n )^{-2}.
\end{equation}
That is, there are $\gg n ( \log n )^{-2}$ elements of $\mathcal A$
that have no prime divisors smaller than $n^{2/9}$. Since the numbers 
in $\mathcal A$ do not exceed $n$, the elements of $\mathcal A$
counted on the left side of \eqref{sieve} have at most four prime
divisors each, that is, the left side of \eqref{sieve} counts solutions
of $n - p = P_4$.

\medskip

We have some freedom in our choice of parameters in the above argument.
For example, we could have set $z = n^{\alpha}$, where $\alpha$ is any
fixed real number in the range $1/5 < \alpha < 1/4$. Of course, what we
would really like to do is set $z = n^{\alpha}$, where $\alpha > 1/4$.
With such a choice for $z$, the above argument would establish the
existence of infinitely many solutions to $n - p = P_3$. Unfortunately,
our choice of $z$ is restricted (via the condition $s = \log y / \log z 
> 2$) by the largest value of $Q$ admissible in the Bombieri--Vinogradov 
theorem. In particular, in order to be able to choose $z = n^{1/4}$,
we would need a version of the Bombieri--Vinogradov theorem that holds
for $Q \le x^{1/2 + \eps}$.

\subsection{Weighted sieves and Chen's theorem}
\label{siv3}

The idea of a \emph{weighted sieve} was introduced by
Kuhn~\cite{Kuhn41} who observed that instead of the sifting
function $S(\mathcal A, \mathfrak P, z)$ one may consider a
more general sum of the type
\begin{equation} \label{iv30}
  W(\mathcal A, \mathfrak P, z) =
  \sum_{ \substack{ a \in \mathcal A,\\ (a, P(z)) = 1}} w(a),
\end{equation}
where $w(a)$ are weights at one's disposal to choose. It is
common to use weights of the form
\begin{equation} \label{iv31}
  w(a) =  1 -
  \sum_{ \substack{ p \mid a \\ z \le p < z_1}} \omega_p,
\end{equation}
with suitably chosen $0 \le \omega_p < 1$. With such a choice of
$w(a)$, \eqref{iv30} can be written in the form
\begin{equation} \label{iv32}
  W(\mathcal A, \mathfrak P, z) = S(\mathcal A, \mathfrak P, z) -
  \sum_{z \le p < z_1} \omega_p S(\mathcal A_p, \mathfrak P, p).
\end{equation}
We can now use an ordinary sieve to estimate the right side of
\eqref{iv32}. For example, we can appeal to \eqref{iv23} to
bound $S(\mathcal A, \mathfrak P, z)$ from below and to
\eqref{iv22} to bound each sifting function $S(\mathcal A_p,
\mathfrak P, p)$ from above. If the resulting lower bound for
the right side of \eqref{iv32} is positive, we then conclude
that there exist elements $a$ of $\mathcal A$ with $w(a) > 0$.
Such numbers $a$ have no prime divisors $p < z$ and the number
of their prime divisors with $z \le p < z_1$ can be controlled
via the choice of the $\omega_p$'s.

The above idea plays an important role in improvements on the
result established in \S \ref{siv2}. Using weighted sieves,
Buchstab~\cite{Buch67} and Richert~\cite{Rich69} proved that
every sufficiently large even $n$ can be represented as the sum
of a prime and a $P_3$-number. Richert used weights of the form
\[
  w(a) =  1 - \theta
  \sum_{ \substack{ p \mid a \\ z \le p < z_1}}
  \left( 1 - \frac {\log p}{\log z_1} \right),
\]
while Buchstab's weights were somewhat more complicated.
Chen's proof of Theorem \ref{th2} uses weights of the form
\eqref{iv31} with $z = n^{1/10}$, $z_1 = n^{1/3}$, and
\[
  \omega_p =  \frac 12 +  \frac 12 \delta_p(a),
\]
where
\[
  \delta_p(a) = \begin{cases}
    1 & \text{if } \; a = pp_1p_2 \; \text{ with } \; p_1 \ge z_1, \\
    0 & \text{otherwise}.
  \end{cases}
\]
Here, $n$ is the even number appearing in the statement of
Theorem \ref{th2} and $a$ is an element of the sequence \eqref{1.7}.
With this choice of $\omega_p$, successful sifting produces
numbers $a \in \mathcal A$ with $w(a) > 0$ and no prime divisors
$p < n^{1/10}$. One can prove that any such number $a$ must 
in fact be a $P_2$-number. The reader can find a detailed
proof of Chen's theorem in \cite[Chapter 11]{HaRi74}, 
\cite[Chapter 10]{Nath96a}, or \cite[Chapter 9]{PaPa92}.

\subsection{Other sieve methods}
\label{siv4}

We conclude our discussion of sieve methods with a brief
account of some of the important ideas in sieve theory
left out of the previous sections.

\paragraph{Selberg's sieve.}

The Rosser--Iwaniec sieve defined by \eqref{iv17} and
\eqref{iv18} is not particularly sensitive to the arithmetical
nature of the sequence $\mathcal A$ that is being sifted. In
fact, the only piece of information about $\mathcal A$ that
the Rosser--Iwaniec sieve does take into account is its sieve
dimension.                                                          %7dec
%$\kappa$ (which is assumed equal to $1$).                          %7dec
Such sieves                                                         %7dec
are known as \emph{combinatorial}. Selberg \cite{Selb47}
proposed another approach, which uses the multiplicative
function $\omega(d)$ appearing in \eqref{iv5} to construct
essentially best possible upper sieve weights $\lambda^+(d)$
for a given sequence $\mathcal A$. 

Suppose that $\rho(d) $ is a real function such that $\rho(1) = 1$. Then
\[
  \sum_{ d \mid n } \mu(d) \le \left( \sum_{ d \mid n } \rho(d) \right)^2.
\]
We can apply this inequality to estimate $S( \mathcal A, \mathfrak P, z )$
as follows:
\begin{align*}
  S( \mathcal A, \mathfrak P, z ) 
  &\le \sum_{ a \in \mathcal A } \left( \sum_{ d \mid n } \rho(d) \right)^2
  = \sum_{ a \in \mathcal A }  \sum_{ d_1, d_2 \mid n } \rho(d_1) \rho(d_2)\\
  &= \sum_{ d_1, d_2} \rho(d_1) \rho(d_2) 
  \big| \mathcal A_{ [d_1, d_2] } \big| ,
\end{align*}
where $ |\mathcal A_d| $ is as before and $[d_1, d_2]$ is the least common
multiple of $d_1$ and $d_2$. Using \eqref{iv5}, we find that
\[
  S( \mathcal A, \mathfrak P, z )  \le XW  + \mathcal R' ,
\]
where
\[
  W =  \sum_{d_1, d_2} \rho(d_1)\rho(d_2)
  \frac{ \omega( [d_1, d_2] ) } { [d_1, d_2 ] } , 
  \qquad      
  \mathcal R' = \sum_{ d_1, d_2} \rho(d_1) \rho(d_2)  r(X, [d_1, d_2] ) .
\] 
In order to control the ``error term'' $\mathcal R'$, we further assume that
$\rho(d) = 0$ when $d > \xi$, where $\xi > 0$ is a parameter. The double sum
$W$ is a quadratic form in the variables $\rho(d)$, $1 < d \le \xi$. Selberg's
idea is to choose the values of these variables as to minimize this quadratic
form.

More information about Selberg's sieve---including the techniques used to 
construct the lower sieve function $\lambda^-(d)$ of Selberg's sieve---can 
be found in \cite{HaRi74, Moto83} and in Selberg's collected works 
\cite{Selb89, Selb91}.

\paragraph{The large sieve.} The method known as the \emph{large
sieve} was introduced in 1941 by Linnik \cite{Linn41}, but its 
systematic study did not commence until R\'enyi's work \cite{Reny47} 
on the binary Goldbach problem. The original idea of Linnik and 
R\'enyi evolved into a general analytic principle that has penetrated 
analytic number theory on many levels (and perhaps does not warrant 
the name ``sieve'' anymore, but the term has survived for historical 
reasons). The most prominent application of the large sieve is 
the Bombieri--Vinogradov theorem. The reader will find discussion 
of the number-theoretic aspects of the large sieve in 
\cite{Bomb74, Dave00, Mont71} and of the analytic side of the 
story in \cite{Dave00, Mont71, Mont78}.

\paragraph{Alternative form of the error term in the sieve.}
Iwaniec~\cite{Iwan80b} obtained a variant of the linear sieve
featuring an error term that is better suited for certain
applications than the error term $\mathcal R$ defined in
\eqref{iv13}. It is of the form
\begin{equation}\label{iv100}
  \sum_{ \substack{ m < M\\ m \mid P(z)}} 
  \sum_{ \substack{ n < N\\ n \mid P(z)}} a_m b_n r(X, mn),
\end{equation}
where the coefficients $a_m$ and $b_n$ are bounded above in 
absolute value and $r(X, mn)$ are the remainder terms defined earlier. 
In some applications, one can use the bilinearity of this expression 
to estimate the double sum when the product $MN$ is larger than the 
largest value of $y$ for which one can obtain a satisfactory bound 
for $\mathcal R$. Iwaniec \cite{Iwan78} 
used this idea in his proof that certain quadratic polynomials take 
on infinitely $P_2$-numbers (recall \S \ref{ssiii4}).

\paragraph{Prime detecting sieves.}

For a long time it was believed that sieve methods are not 
capable of detecting prime numbers; there are even a couple of 
prominent papers (see \cite{Bomb76, Selb89}) that quantify the
shortcomings of the classical sieve technology. In short, classical
sieves are incapable of distinguishing between integers having even
number of prime divisors and those having an odd number of prime
divisors (this is known in sieve theory as the \emph{parity 
obstacle}). A prime detecting sieve overcomes the parity obstacle
by combining the general sieve philosophy with additional analytic
information. A variant of the basic idea can be traced all the 
way back to Vinogradov's work on sums over primes, but the first 
explicit uses of prime detecting sieves appeared in the late 
1970s in investigations of the distribution of primes in short 
intervals (see \cite{H-BIw79, IwJu79}). The method flourished
during the last decade and has been instrumental in the proofs
of several of the restults mentioned in the previous sections:
the result of Friedlander and Iwaniec \cite{FrIw98a} on prime
values of $x^2 + y^4$; the results of Heath-Brown and Moroz
\cite{H-Br01, H-BMo02} on prime values of binary cubic forms;
and the result of Baker, Harman, and Pintz \cite{BaHaPi01} on 
primes in short intervals are just three such examples. Compared
to classical sieve methods, the theory of prime detecting sieves 
is still in its infancy and thus the general literature on the
subject is relatively scarce, but the reader eager to learn more
about such matters will find two excellent expositions in
\cite{FrIw98b} and \cite{Harm97}.

\section{Other work on the Waring--Goldbach problem}
\label{sv}

In the Introduction, we mentioned the cornerstones in the study of
the Goldbach and Waring--Goldbach problems. However, as is often the
case in mathematics, those results are intertwined with a myriad of
other results on various aspects and variants of the two main 
problems. In this final section, we describe some of the more
important results of the latter kind. The circle method, sieve 
methods, or a combination of them play an essential role in the 
proofs of all these.

\subsection{Estimates for exceptional sets}
\label{ES}

Inspired by the work of Chudakov~\cite{Chud38}
and Estermann~\cite{Este38} on the exceptional
set in the binary Goldbach problem, Hua studied the function 
$h(k)$, defined to be the least $s$ such that almost all integers
$n \le x$, $n \equiv s \pmodulo {K(k)}$, can be written as the
sum of $s$ $k$th powers of primes
($K(k)$ is defined by \eqref{1.12}). Let $E_{k, s}(x)$ denote the
number of exceptions, that is, the number of integers $n$, with 
$n \le x$ and $n \equiv s \pmodulo {K(k)}$, for which \eqref{1.14}
has no solution in primes $p_1, \dots, p_s$. Hua showed
(essentially) that if $H(k) \le s_0(k)$, then $E_{k, s}(x) = o(x)$
for any $s \ge \frac 12 s_0(k)$. Later, Schwarz \cite{Schw61}
refined Hua's method to show that
\begin{equation}\label{i16}
  E_{k, s}(x) \ll x(\log x)^{-A}
\end{equation}
for any fixed $A > 0$.

In recent years, motivated by the estimate \eqref{1.6} of Montgomery
and Vaughan, several authors have pursued similar estimates for 
exceptional sets for squares and higher powers of primes. The first
to obtain such an estimate were Leung and Liu~\cite{LeLi93}, who
showed that $E_{2, 3}(x) \ll x^{1 - \delta}$, with
an absolute constant $\delta > 0$. Explicit versions of this result
were later given in~\cite{BaLiZh00, HaKu04, Kumc02b, LiZh01, LiZh03},
the best result to date being the estimate (see Harman and Kumchev
\cite{HaKu04})
\[
  \textstyle E_{2, 3}(x) \ll x^{6/7 + \eps}.
\]
Furthermore, several authors~\cite{HaKu04, JLiu03, LiLi00b,
LiWoYu03, Wool02b} obtained improvements on Hua's bound
\eqref{i16} for $E_{2, 4}(x)$, the most recent being the bound
\[
  \textstyle E_{2, 4}(x) \ll x^{5/14 + \eps},
\]
established by Harman and Kumchev \cite{HaKu04}.
Ren~\cite{Ren00} studied the exceptional set for sums of
five cubes of primes and proved that
\[
  \textstyle E_{3, s}(x) \ll x^{1 - (s - 4)/153}
  \qquad (5 \le s \le 8). 
\]
This estimate has since been improved by Wooley~\cite{Wool02a}
and Kumchev~\cite{Kumc02a}. In particular, Kumchev \cite{Kumc02a}
showed that
\[
  \begin{split}
    \textstyle
    E_{3, 5}(x) \ll x^{79/84}, \qquad E_{3, 6}(x) \ll x^{31/35}, \\
    \textstyle
    E_{3, 7}(x) \ll x^{51/84}, \qquad E_{3, 8}(x) \ll x^{23/84}.
  \end{split}
\]
Finally, Kumchev~\cite{Kumc02a} has developed the
necessary machinery to obtain estimates of the form $E_{k, s}(x)
\ll x^{1 - \delta}$, with explicit values of $\delta =
\delta(k, s) > 0$, for all pairs of integers $k \ge 4$ and $s$
for which an estimate of the form \eqref{i16} is known.

In 1973 Ramachandra~\cite{Ramach73} considered the exceptional
set for the binary Goldbach problem in short intervals.
He proved that if $y \ge x^{7/12 + \eps}$ and $A > 0$, then
\[
  E(x + y) - E(x) \ll y(\log x)^{-A},
\]
where the implied constant depends only on $A$ and $\eps$.
After a series of improvements on this result~\cite{BaHaPi97,
Duf94, Duf95, Jia95a, Jia95a+, Jia96a, HLi95, Mika92b, PerPin93},
this estimate is now known for $y \ge x^{7/108 + \eps}$ (see
Jia~\cite{Jia96a}). Lou and Yao~\cite{LouYao81, Yao82} were the
first to pursue a short interval version of the estimate
\eqref{1.6} of Montgomery and Vaughan. 
Their result was substantially improved by Peneva \cite{Pene01}
and the best result in this
direction, due to Languasco~\cite{Langu04}, states that
there exists a small constant $\delta > 0 $ such that
\[
  E(x + y) - E(x) \ll y^{1 - \delta / 600},
\]
whenever $y \ge x^{7/24 + 7 \delta}$. 

Furthermore, J. Liu and Zhan~\cite{LiZh97b} and 
Mikawa~\cite{Mika97} studied the quantity $E_{2, 3}(x)$ in 
short intervals and the latter author showed that
\[
  E_{2, 3}(x + y) - E_{2, 3}(x) \ll y(\log x)^{-A}
\]
for any fixed $A>0$ and any $y \ge x^{1/2 + \eps}$.

\subsection{The Waring--Goldbach problem with almost \\ primes}
\label{QWGp}

There have also been attempts to gain further knowledge about
the Waring--Goldbach problem by studying closely related but
more accessible problems. The most common such variants relax
the multiplicative constraint on (some of) the variables.
Consider, for example, Lagrange's equation
\begin{equation}\label{i17}
  x_1^2 + x_2^2 + x_3^2 + x_4^2 = n.
\end{equation}
Greaves~\cite{Grea76} proved that every sufficiently large $n \not
\equiv 0, 1, 5 \pmodulo 8$ can be represented in the form
\eqref{i17} with $x_1, x_2$ primes and $x_3, x_4$ (unrestricted)
integers. Later, Plaksin~\cite{Plak82} and Shields~\cite{Shie79}
found independently an asymptotic formula for the number of such
representations. Br\"udern and Fouvry~\cite{BrFo94} proved that
every sufficiently large integer $n \equiv 4 \pmodulo {24}$ can be
written as the sum of four squares of $P_{34}$-numbers. Heath-Brown
and Tolev~\cite{H-BTo03} established, under the same hypothesis on
$n$, that one can solve \eqref{i17} in one prime and three almost
primes of type $P_{101}$ or in four almost primes, each of type
$P_{25}$. Tolev~\cite{Tole03} has recently improved the results in
\cite{H-BTo03}, replacing the types of the almost primes involved
by $P_{80}$ and $P_{21}$, respectively. We must also mention
the recent result by Blomer and Br\"udern \cite{BlBr04} that all
sufficiently large integers $n$ such that $n \equiv 3 \pmodulo{24}$
and $5 \nmid n$ are sums of three almost primes
of type $P_{521}$ (and of type $P_{371}$ if $n$ is also squarefree).

In 1951 Roth \cite{Roth51} proved that if $n$ is sufficiently
large, the equation
\begin{equation}\label{5.3}
  x^3 + p_1^3 + \dots + p_7^3 = n
\end{equation}
has solutions in primes $p_1, \dots, p_7$ and an integer $x$.
Br\"udern~\cite{Brud95a} showed that if $n \equiv 4 \pmodulo{18}$,
then $x$ can be taken to be a $P_4$-number, and
Kawada \cite{Kawa97} used an idea from Chen's proof of Theorem
\ref{th2} to obtain a variant of Br\"udern's result for almost
primes of type $P_3$. Furthemore, Br\"udern~\cite{Brud95b}
proved that every sufficiently large integer is the sum of the
cubes of a prime and six almost-primes (five $P_5$-numbers and a
$P_{69}$-number) and Kawada \cite{Kawa??} has shown that every
sufficiently large integer is the sum of seven cubes of $P_4$-numbers.   %7dec

Wooley~\cite{Wool02b} showed that all but $O\big( (\log x)^{6 +
\eps} \big)$ integers $n \le x$, satisfying certain natural
congruence conditions can be represented in the form \eqref{i17}
with prime variables $x_1, x_2, x_3$ and an integer $x_4$.
Tolev~\cite{Tole03+} established a result of similar strength
for the exceptional set for equation \eqref{i17} with primes
$x_1, x_2, x_3$ and an almost prime $x_4$
of type $P_{11}$.

\subsection{The Waring--Goldbach problem with restricted \\ variables}
\label{WGprv}

Through the years, a number of authors have studied variants of
the Goldbach and Waring--Goldbach problems with additional
restrictions on the variables. In 1951 Haselgrove \cite{Hase51}
announced that every sufficiently large odd integer $n$ is the
sum of three primes $p_1, p_2, p_3$ such that $|p_i - n/3| \le
n^{63/64 + \eps}$. In other words, one can take the primes in
Vinogradov's three prime theorem to be ``almost equal''.
Subsequent work by several mathematicians \cite{BaHa98, Chen65, 
Jia89, Jia94, PaPa89, Zhan91} tightened the range for the $p_i$'s to 
$|p_i - n/3| \le n^{4/7}$ (see Baker and Harman~\cite{BaHa98}).

Furthermore, Bauer, Liu, and Zhan~\cite{Bauer97, LiZh96,
LiZh00} considered the problem of representations of an
integer as sums of five squares of almost equal primes.
The best result to date is due to Liu and Zhan \cite{LiZh00},
who proved that every sufficiently large integer $n \equiv 5
\pmodulo {24}$ can be written as
\[
  n = p_1^2 + \dots + p_5^2,
\]
with primes $p_1, \dots, p_5$ satisfying $| p_i^2 - n/5 |
< n^{45/46 + \eps}$. Liu and Zhan \cite{LiZh96} also showed
that the exponent $\frac {45}{46}$ can be replaced by
$\frac {19}{20}$ on the assumption of GRH.

In 1986 Wirsing~\cite{Wirs86} proved that there exist sparse
sequences of primes $\mathcal S$ such that every sufficiently
large odd integer can be represented as the sum of three
primes from $\mathcal S$. However, his method was
probabilistic and did not yield an example of such a sequence.
Thus, Wirsing proposed the problem of finding ``natural'' 
examples of arithmetic sequences having this
property. The first explicit example was given by Balog and
Friedlander~\cite{BaFr92}. They proved that the sequence
of Piatetski-Shapiro primes (recall \eqref{thinset2})
is admissible for $1 < c < 21/20$. Jia~\cite{Jia95b}
improved the range for $c$ to $1 < c < 16/15$, and
Peneva~\cite{Pene03} studied the binary problem with a
Piatetski-Shapiro prime and an almost prime.
Tolev~\cite{Tole99}--\cite{Tole00+} and
Peneva~\cite{Pene00} considered additive problems with
prime variables $p$ such that the integers $p + 2$ are
almost-primes. For example, Tolev \cite{Tole00+} proved
that every sufficiently large $n \equiv 3 \pmodulo {6}$
can be represented as the sum of primes $p_1, p_2, p_3$
such that $p_1+2=P_2$, $p_1+2=P_5$, and $p_1+2=P_7$.
Green and Tao announced at the end of \cite{GrTa04} that, 
using their method, one can prove that there are arbitrarily long 
non trivial arithmetic progressions consisting of primes $p$
such that $p+2 = P_2$. They
presented in \cite{GrTa04+} a proof of this result
for progressions of three primes.

\subsection{Linnik's problem and variants}
\label{Lpv}

In the early 1950s Linnik proposed the problem of finding
sparse sequences $\mathcal A$ such that all sufficiently
large integers $n$ (possibly subject to some parity condition)
can be represented as sums of two primes and an element of
$\mathcal A$. He considered two special sequences. First,
he showed \cite{Linn52} that if GRH holds, then every
sufficiently large odd $n$ is the sum of three primes $p_1,
p_2, p_3$ with $p_1 \ll (\log n)^3$. Montgomery and Vaughan
\cite{MoVa75} sharpened the bound on $p_1$ to $p_1 \ll
(\log n)^2$ and also obtained an unconditional result with
$p_1 \ll n^{7/72 + \eps}$; the latter bound has been
subsequently improved to $p_1 \ll n^{0.02625}$ (this follows
by the original argument of Montgomery and Vaughan from
recent results of Baker, Harman, and Pintz \cite{BaHaPi01}
and Jia \cite{Jia96b}).

Linnik \cite{Linn51, Linn53} was also the first to study additive
representations as sums of two primes and a fixed number of powers 
of $2$. He proved, first under GRH and later unconditionally, that 
there is an absolute constant $r$ such that every sufficiently 
large even integer $n$ can be expressed as the sum of two primes 
and $r$ powers of $2$, that is, the equation
\[
  p_1 + p_2 + 2^{\nu_1} + \dots + 2^{\nu_r} = n,
\]
has solutions in primes $p_1, p_2$ and non-negative integers 
$\nu_1, \dots, \nu_r$. Later Gallagher~\cite{Gala75} established 
the same result by a different method. Several authors have used
Gallagher's approach to find explicit values of the constant $r$
above (see \cite{HLi01, HLi01a, LiLiWa98a, LiLiWa98b, LiLiWa99,
Wang99}); in particular, Li~\cite{HLi01a} proved that $r = 1906$ 
is admissible and Wang \cite{Wang99} obtained $r = 160$ under GRH. 
Recently, Heath-Brown and Puchta~\cite{H-BPu02} and Pintz and 
Ruzsa~\cite{PiRu03} made (independently) an important discovery 
that leads to a substantial improvement on the earlier results. 
Their device establishes Linnik's result with $r = 13$ (see 
\cite{H-BPu02}) and with $r = 7$ under GRH (see \cite{H-BPu02,
PiRu03}). Furthermore, Pintz and Rusza \cite{PiRu04} have announced 
an unconditional proof of the case $r = 8$.

There is a similar approximation to the Waring--Goldbach problem
for four squares of primes. J. Y. Liu, M. C. Liu, and Zhan 
\cite{LiLiZh99, LiLiZh02} proved that there exists a constant $r$ 
such that every sufficiently large even integer $n$ can be expressed 
in the form
\[
  p_1^2 + p_2^2 + p_3^2 + p_4^2 + 2^{\nu_1} + \dots + 2^{\nu_r} = n,
\]
where $p_1, \dots, p_4$ are primes and $\nu_1, \dots, \nu_r$ are 
non-negative integers. J. Y. Liu and M. C. Liu \cite{LiLi00a} 
established this result with $r = 8330$ and considered also the 
related problem about representations of integers as sums of a prime, 
two squares of primes and several powers of 2.

\subsection{Additive problems with mixed powers}
\label{OAP}

In 1923 Hardy and Littlewood~\cite{HaLi23b} used the general
philosophy underlying the circle method to formulate several 
interesting conjectures. For example, they stated a conjectural 
asymptotic formula for the number of representations of a large 
integer $n$ in the form
\begin{equation}\label{d4}                                          
  p + x^2 + y^2 = n,                                                
\end{equation}                                                       
where $p$ is a prime and $x, y$ are integers. Their prediction
was confirmed in the late 1950s, first by Hooley~\cite{Hoo57} 
under the assumption of GRH and then unconditionally by Linnik
\cite{Linn60}. The reader will find the details of the proof in 
\cite{Hoo76, Linn63}.

In another conjecture, Hardy and Littlewood proposed an asymptotic
formula for the number of representations of a large integer $n$
as the sum of a prime and a square. While such a result appears to
lie beyond the reach of present methods, Miech \cite{Mie68b} showed 
that this conjecture holds for almost all integers $n \le x$. Let 
$E_k(x)$, $k \ge 2$, denote the number of integers $n \le x$ such 
that the equation $n = p + x^k$ has no solution in a prime $p$ and 
an integer $x$. Miech obtained the bound $E_2(x) \ll x(\log x)^{-A}$ 
for any fixed $A > 0$. Subsequent work of Br\"udern, Br\"unner, 
Languasco, 
Mikawa, Perelli, Pintz, Polyakov, A. I. Vinogradov, and Zaccagnini
\cite{BrPer96, BrunPerPi89, Langu04+,
Mika93, PerZac95, Pol90, AVin85, Zacc92} 
extended and sharpened Miech's estimate considerably. Here is a 
list of some of their results:
\begin{itemize}
  \item For any fixed $k \ge 2$, we have 
    $E_k(x) \ll x^{1 - \delta_k}$, where $\delta_k > 0$ depends
    at most on $k$; see \cite{BrunPerPi89, Pol90, AVin85} for 
    the case $k = 2$ and \cite{PerZac95, Zacc92} for the general 
    case.
  \item Assuming GRH, we have $E_k(x) \ll x^{1 - \delta_k}$,
    where $\delta_k = 1 / (k2^k)$ or $\delta_k = 1 / (25k)$
    according as $2 \le k \le 4$ or $k \ge 5$; see 
    \cite{PerZac95} and \cite{BrPer96}.
  \item If $ k \ge 2 $ is a fixed integer and $ K = 2^{k-2} $
   then there exists a small absolute constant $ \delta > 0 $
   such that
    \[
      E_k(x + y) - E_k(x) \ll y^{ 1 - \delta / (5K) } ,
    \]
    provided that $x^{(7/12)(1 - 1/k) + \delta} \le y \le x$;
    see \cite{Langu04+}.
\end{itemize}
Furthermore, several mathematicians \cite{Bauer98, Bauer99, 
BrPer96, LiZh97b} obtained variants of the above bounds
in the case when the variable $x$ is also restricted to primes,
while Zaccagnini \cite{Zacc01} studied the more general problem
of representing a large integer $n$ in the form $n = p + f(x)$,
where $f(X) \in \mathbb Z[X]$.

Several interesting theorems were proved by
Br\"udern and Kawada~\cite{BrKa02}.
For example, one of them states that if $k$
is an integer with $3 \le k \le 5$, then all
sufficiently large integers $n$ can be represented as
\[
   x + p_1^2 + p_2^3 + p_3^k = n ,
\]
where $p_i$ are primes and $x=P_2$.

\subsection{The Waring--Goldbach problem ``with coefficients''}
\label{s5.6}

In this section we discuss the solubility of equations of the form \eqref{1.16}, 
which we introduced in \S\ref{s1.4} as a natural generalization of the
Waring--Goldbach 
problem. There are two substantially different contexts in which one can study this
problem. 
Suppose first that all $a_1, \dots, a_s, n$ are all of the same sign. 
Then one expects that \eqref{1.16} must have solutions for sufficiently large $|n|$. 
When ``sufficiently large'' is understood as $|n| \ge C(a_1, \dots, a_s)$, with some 
unspecified constant depending on the $a_j$'s, this is a trivial modification of the 
Waring--Goldbach problem (that can be handled using essentially the same tools). 
On the other hand, the problem of finding solution when $|n|$ is not too large 
compared to $|\mathbf a|_{\infty} = \max\{ |a_1|, \dots, |a_s| \}$ is significantly 
more challanging. Similarly, if $a_1, \dots, a_s$ are not all of the same sign, 
one wants to find solutions of \eqref{1.16} in primes $p_1, \dots, p_s$ that are 
not too large compared to $|\mathbf a|_{\infty}$ and $|n|$. Such questions were 
investigated first by Baker \cite{Bake67}, who studied the case $k = 1$ and $s = 3$. 
Later, Liu and Tsang \cite{LiTs89} showed, again for $k = 1$ and $s = 3$, that 
\eqref{1.16} has solutions when:
\begin{itemize}
  \item $a_1, a_2, a_3$ are of the same sign and $|n| \gg |\mathbf a|_{\infty}^A$ 
  for some absolute constant $A > 0$;
  \item $a_1, a_2, a_3$ are not of the same sign and 
  $\max \{ p_1, p_2, p_3 \} \ll |\mathbf a|_{\infty}^{A - 1} + |n|$.
\end{itemize}
In these results, the coefficients $a_1, a_2, a_3, n$ must satisfy also certain
necessary congruence conditions (which generalize the requirement that $n$ be
odd in Vinogradov's three primes theorem). Through the efforts of several 
mathematicians, the constant $A$ has been evaluated and it is known that the 
value $A = 38$ is admissible (see Li \cite{HLi01b}). Furthermore, if we 
replace the natural arithmetic conditions on the coefficients by another set
of conditions, which are somewhat more restrictive but also simplify greatly
the analysis, we can decrease the value of $A$ further. In particular, Choi 
and Kumchev \cite{ChKuX1} have shown that $A = 23/3$ is admissible under 
such stronger hypotheses.

Liu and Tsang \cite{LiTs91} studied also the quadratic case of \eqref{1.16} in five
variables and obtained results similar to those stated above for the linear case. In
this problem, explicit values of the analogue of $A$ above were given by Choi and
Liu~\cite{ChLi02, ChLiX1}, Choi and Kumchev \cite{ChKu04}, and Harman and Kumchev
\cite{HaKu04}.
In particular, it is proved in \cite{HaKu04} that \eqref{1.16} with $k = 2$ and $s =
5$ has solutions when:
\begin{itemize}
  \item $a_1, \dots, a_5$ are of the same sign and $|n| \gg |\mathbf a|_{\infty}^{15
+ \eps}$;
  \item $a_1, \dots, a_5$ are not of the same sign and $\max\{ p_1, \dots, p_5 \}
\ll |\mathbf a|_{\infty}^{7 + \eps} + |n|^{1/2}$.
\end{itemize}

\subsection{Diophantine inequalities with primes}
\label{DI}

Some variants of the Waring--Goldbach problem are stated most
naturally in terms of diophantine inequalities. The best-known 
problem of this kind concerns the distribution of the values 
of the forms
\begin{equation}\label{DI1}
  \lambda_1 p_1^k + \dots + \lambda_s p_s^k,
\end{equation}
where $k$ and $s$ are positive integers, $\lambda_1, \dots,
\lambda_s$ are nonzero real numbers, and $p_1, \dots, p_s$ are
prime variables. It is natural to conjecture that if $\lambda_1, 
\dots, \lambda_s$ are not all of the same sign and if $\lambda_i 
/ \lambda_j$ is irrational for some pair of indices $i, j$, then 
the values attained by the form \eqref{DI1} are dense in 
$\mathbb R$ whenever $s \ge s_0(k)$. In other words, given any 
$\eps > 0$ and $\alpha \in \mathbb R$, the inequality
\begin{equation}\label{DI2}
  \left| \lambda_1 p_1^k + \dots + \lambda_s p_s^k - \alpha \right|
  < \eps
\end{equation}
should have a solution in primes $p_1, \dots, p_s$. The first 
results in this problem were obtained by Schwarz \cite{Schw63},
who established the solvability of \eqref{DI2} under the same
restrictions on $s$ as in Theorem \ref{th3}. Baker \cite{Bake67}
and Vaughan \cite{Vaug74a, Vaug74b, Vaug76} proposed the more
difficult problem of replacing the fixed number $\eps$ on the right 
side of \eqref{DI2} by an explicit function of $\max\{p_1, \dots, 
p_s\}$ that approaches~$0$ as $\max\{p_1, \dots, p_s\} \to \infty$. 
Further work has focused primarily on the case of small $k$. For 
example, Harman \cite{Harm90} has shown that under the above 
assumptions on $\lambda_1, \lambda_2, \lambda_3$, the diophantine 
inequality
\[
  \big| \lambda_1 p_1 + \lambda_2 p_2 + \lambda_3 p_3 - 
  \alpha \big| < \max\{p_1, p_2, p_3\}^{-1/5 + \eps}
\]
has infinitely many solutions in primes $p_1, p_2, p_3$. Baker and 
Harman \cite{BaHa82} showed that on GRH the exponent $\frac 15$ in
this result can be replaced by $\frac 14$. Furthermore, Harman 
\cite{Harm84} proved that if $\lambda_1/\lambda_2$ is a negative
irrational number, then for any real $\alpha$ the inequality
\[
  \big| \lambda_1 p + \lambda_2 P_3 - \alpha \big| < p^{ -1/300 }
\]
has infinitely many solutions in a prime $p$ and a $P_3$-almost prime.
(This improves on an earlier result of Vaughan \cite{Vaug76}, where the
almost prime is a $P_4$-number.) 

In 1952 Piatetski-Shapiro~\cite{PS52} considered a variant of
the Waring--Gold\-bach problem for non-integer exponents $c > 1$.
He showed that for any fixed $c > 1$, which is not an integer, 
there exists an integer $H(c)$ with the following property: if
$s \ge H(c)$,  the inequality           
\begin{equation} \label{PS}                                         
  \big| p_1^c + \dots + p_s^c - \alpha \big| < \eps                         
\end{equation}                                                      
has solutions in primes $p_1, \dots, p_s$ for any fixed $\eps > 0$ 
and $\alpha \ge \alpha_0(\eps, c)$. In particular, Piatetski-Shapiro
showed that $H(c) \le 5$ for $1 < c < 3/2$. Motivated
by Vinogradov's three prime theorem, Tolev \cite{Tole92} proved that 
$H(c) \le 3$ for $1 < c < 15/14$. The range of validity of Tolev's
result was subsequently extended by several authors \cite{Cai96a, 
Cai96b, Kumc99, KumcNede98}; in particular, Kumchev \cite{Kumc99} 
has given the range $1 < c < 61/55$. Furthemore, it follows from 
the work of Kumchev and Laporta~\cite{KumcLapo02, Lapo99c} that 
$H(c) \le 4$ for $1 < c < 6/5$ and for almost all (in the sense of 
Lebesgue measure) $1 < c < 2$, while Garaev~\cite{Gara03} has 
showed that $H(c) \le 5$ for $1 < c < (1 + \sqrt 5)/2 = 1.61\dots$.
Finally, Tolev \cite{Tole95} and Zhai~\cite{Zhai00} have studied
systems of inequalities of the form \eqref{PS}.
                           
Several authors \cite{ArBuCh99, ArCheCh03, Bu87, Bu89} have
studied variants of Goldbach's problem, suggested by results 
about additive inequalities. For example, Arkhipov, Chen, and 
Chubarikov~\cite{ArCheCh03} proved that if $\lambda_1/\lambda_2$ 
is an algebraic irrationality, then all but $O(x^{2/3 +\eps})$ 
positive integers $n\le x$ can be represented in the form 
\[
 [\lambda_1 p_1] + [\lambda_2 p_2] = n,
\]
where $p_1, p_2$ are primes.

\section{A new path: arithmetic progressions of \\ primes} 

Finally, we should say a few words about the astonishing result 
of Green and Tao~\cite{GrTa04} on the existence of arbitrarily
long arithmetic progressions of prime numbers. They deduce the 
existence of such arithmetic progressions from a generalization 
of a celebrated theorem of Szemer\'edi  \cite{Sze69, Sze75}, 
which is itself a deep result in combinatorial number theory. 
Let $\mathcal A$ be a set of positive integers with 
\emph{positive upper density}, that is,
\[
  \delta(\mathcal A) = \limsup_{N \to \infty} 
  \frac {\#\{ n \in \mathcal A : n \le N\}}N > 0.
\]
In its original, most basic form, Szemer\'edi's theorem asserts that 
such a set $\mathcal A$ contains an arithmetic progression of length 
$k$ for all integers $k \ge 3$. From this basic statement, Green and
Tao deduce the following more general result.

\begin{theorem}[Szemer\'edi's theorem for pseudorandom measures]\label{Sze}
  Let $\delta \in (0, 1] $ be a fixed real number, let $k \ge 3 $ be 
  a fixed integer, and let $N$ be a large prime.
  Suppose that $\nu$ is a ``$k$-pseudorandom measure\footnote{A 
  $k$-pseudorandom measure on $\mathbb Z_N$ is a non-negative function
  on $\mathbb Z_N$ whose average over $\mathbb Z_N$ is close to $1$ and
  which is subject to a couple of additional constraints that are too
  technical to state here. See \cite{GrTa04} for details.}'' on 
  $\mathbb Z_N = (\mathbb Z / N\mathbb Z)$ and $f: \mathbb Z_N \to 
  [0, \infty)$ is a function satisfying
  \begin{equation} \label{GT1}
    0 \le f(x) \le \nu (x) \qquad
    \text{for all } \, x \in \mathbb Z_N 
  \end{equation}
  and
  \[
    \sum_{ x \in \mathbb Z_N } f(x) \ge \delta N.
  \]
  Then 
  \begin{equation} \label{GT2}
    \sum_{ x \in \mathbb Z_N } \sum_{ r\in \mathbb Z_N }
    f(x) f(x + r) \dots f( x + (k-1) r ) \gg N^2,
  \end{equation}
  the implied constant depending at most on $\delta$ and $k$.
\end{theorem}

To relate this result to the version
of Szemer\'edi's theorem stated earlier, consider the case where
$\nu(x) = 1$ for all $x$ (this is a $k$-pseudorandom measure) 
and $f(x)$ is the characteristic function of the set
$\mathcal A_N = \mathcal A \cap [1, N]$ considered as a subset 
of $\mathbb Z_N$. Then the left side of \eqref{GT2} counts 
(essentially) the $k$-term arithmetic progressions in the set 
$\mathcal A_N$ (the majority of which are also $k$-term 
arithmetic progressions in $\mathcal A \cap \mathbb Z$). 

To derive the result on arithmetic progressions of primes, Green 
and Tao take $f(x)$ to be                                                 %7des
a function which,                                                         %7des
in some sense (see \cite{GrTa04} for details), approximates               %7des
%a majorant for a weighted version of                                     %7des 
%                                                                         %7des
%     MISLJA, CHE TAKA KAKTO SYM GO NAPRAVIL                              %7des
%     E MALKO PO-TOCHNO.                                                  %7des  
%                                                                         %7des  
the characteristic function of the primes in the interval 
$[c_1N, c_2N]$, where $0 < c_1 <c_2 <1$ are suitable constants.
Then they construct a pseudorandom measure $\nu(x)$ such that
\eqref{GT1} holds. This leads to the following theorem.

\begin{theorem}[Green and Tao, 2004]
  Let $k \ge 3$ and let $\mathcal A$ be a set of prime numbers 
  such that
  \[
    \limsup_{ N \to \infty } 
    \frac {\#\{ n \in \mathcal A : n \le N\}}{\pi(N)} > 0.
  \]
  Then $\mathcal A$ contains infinitely many $k$-term 
  arithmetic progresions. In particular, there are infinitely many 
  $k$-term arithmetic progresions of prime numbers.
\end{theorem}

We remark that the infinitude of the $k$-term progressions of primes
is a consequence of \eqref{GT2}. In fact, using the explicit form
of the function $f(x)$ to which they apply Theorem \ref{Sze}, Green
and Tao establish the existence of $\gg N^2(\log N)^{-k}$ $k$-term
progressions within $\mathcal A \cap [1, N]$.

Several other interesing results are announced in \cite{GrTa04}.
For example, one of them asserts that there are infinitely many 
progressions of primes $p_1, \dots, p_k$ such that each $p_i + 2 $ 
is a $P_2$-number (a proof of this result in the case $k = 3$ is 
presented in \cite{GrTa04+}).

\begin{conclusion}
  With this, our survey comes to a close. We tried to describe the central problems 
and the main directions of research in the additive theory of prime numbers and 
to introduce the reader to the classical methods. Complete success in such an 
undertaking is perhaps an impossibility, but hopefully we have been able to paint 
a representative picture of the current state of the subject and to motivate the
reader 
to seek more information from the literature. Maybe some of our readers will one day 
join the ranks of the number theorists trying to turn the great conjectures mentioned 
above into beautiful theorems!
\end{conclusion}

\begin{acknowledgement}
  This paper was written while the first author enjoyed the benefits of 
postdoctoral positions at the University of Toronto and at the University 
of Texas at Austin. He would like to take this opportunity to express his 
%deepest                                                                  %7dec
%                                                                         %7dec
%  TUK KOMENTIRAH "DEEPEST" PONEZHE PO-KYSNO SPOMENAVANE                  %7dec
%  IMENATA NA FRIEDLANDER I HEATH-BROWN.                                  %7dec 
%  MOZHE BI E PO-UMESTNO TAZI DUMA DA SE UPOTREBI TAM,                    %7dec
%  A NE ZA INSTITUZIITE. NO AKO NE TI HARESVA -                           %7dec
%  NAPRAVI GO KAKTO TI ISKASH.                                            %7dec 
%  SYSHTO TAKA PROMENIH LEKO TEKSTA ZA BLAGODARNOSTTA,                    %7dec
%  KOITO TI PREDLOZHI.                                                    %7dec
%                                                                         %7dec
gratitude to these institutions for their support. The second author 
was supported by Plovdiv University Scientific Fund grant 03-MM-35.
Last but not least, the authors would like to thank                       %7dec
Professors J. Friedlander                                                 %7dec
and D.R. Heath-Brown for several valuable discussions                     %7dec
over the years and for some useful comments about the preliminary         %7dec
version of the survey.                                                    %7dec
\end{acknowledgement}

\bigskip

\textsc{Department of Mathematics, 
        1 University Station, C1200,
        The University of Texas at Austin,
        Austin, TX 78712,
        U.S.A.}
         
\textit{E-mail:} \texttt{kumchev@math.utexas.edu}

\medskip

\textsc{Department of Mathematics,
        Plovdiv University ``P. Hilendarski'',
        24 Tsar Asen Street,
        Plovdiv 4000,
        Bulgaria}
         
\textit{E-mail:} \texttt{dtolev@pu.acad.bg }

\end{document}